\documentclass{article}
\usepackage{latexsym,amsmath,amssymb,amscd}
\usepackage[dvips]{color}
\numberwithin{equation}{section}

\begin{document}

\title{On slowly rotating axisymmetric solutions of the Einstein-Euler equations}
\author{Tetu Makino \footnote{Professor Emeritus at Yamaguchi University, Japan. E-mail:
makino@yamaguchi-u.ac.jp}}
\date{\today}
\maketitle

\newtheorem{Lemma}{Lemma}
\newtheorem{Proposition}{Proposition}
\newtheorem{Theorem}{Theorem}
\newtheorem{Definition}{Definition}
\newtheorem{Remark}{Remark}
\newtheorem{Corollary}{Corollary}

\begin{abstract}
In recent works we have constructed axisymmetric solutions to the Euler-Poisson equations which give mathematical models of slowly uniformly rotating gaseous stars. We try to extend this result to the study of solutions of the Einstein-Euler equations in the framework of the general theory of relativity. Although many interesting studies have been done about axisymmetric metric in the general theory of relativity, they are restricted to the region of the vacuum. Mathematically rigorous existence theorem of the axisymmetric interior solutions of the stationary metric corresponding to the energy-momentum tensor of the perfect fluid with non-zero pressure may be not yet established until now
except only one found in the pioneering work by U. Heilig done in 1993. In this article, along a different approach to that of  Heilig's work, axisymmetric stationary solutions of the Einstein-Euler equations are constructed near those of the Euler-Poisson equations when the speed of light is sufficiently large in the considered system of units, or, when the gravitational field is sufficiently weak.\\


\end{abstract}

\section{Introduction}

We consider the Einstein's equations
\begin{equation}
R_{\mu\nu}-\frac{1}{2}g_{\mu\nu}R=\frac{8\pi\mathsf{G}}{\mathsf{c}^4}T_{\mu\nu}
\end{equation}
for the metric
\begin{equation}
ds^2=g_{\mu\nu}dx^{\mu}dx^{\nu}
\end{equation}
with the energy-momentum tensor of the perfect fluid
\begin{equation}
T^{\mu\nu}=(\mathsf{c}^2\rho+P)U^{\mu}U^{\nu}-Pg^{\mu\nu}.
\end{equation}
Here $\mathsf{G}$ is the gravitational constant, and $\mathsf{c}$ is the speed of light. They are positive constants. $R_{\mu\nu}$ is the Ricci tensor and $R=g^{\mu\nu}R_{\mu\nu}$ is the scalar curvature. 
The Greek letters $\mu, \nu, ...$ for index take values $0,1,2,3$, and the Einstein's rule on summations is adopted.
For the definitions of these symbols and the meaning of  the fundamental equations of
the gravitational field in the general theory of relativity, see e.g. 
\cite{LandauL}. 

In this article we assume that the pressure $P$ is a given function of
the density $\rho$. More precisely, we put the following assumption:\\

{\bf (A): \   $P$ is a given smooth function of $\rho>0$
such that $0<P, 0<dP/d\rho<\mathsf{c}^2$ for $\rho >0$, and there is a smooth function
$\Lambda \in C^{\infty}(\mathbb{R})$ such that $\Lambda(0)=0$, 
$\Lambda$ is analytic at $0$, and
\begin{equation}
P=\mathsf{A}\rho^{\gamma}(1+\Lambda(\mathsf{A}\rho^{\gamma-1}/\mathsf{c}^2))\quad
\mbox{for}\quad \rho >0.
\end{equation}
Here $\mathsf{A}, \gamma$ are positive constants and 
\begin{equation}
\frac{6}{5}<\gamma < 2.
\end{equation}
}\\

In the works \cite{JMKU98}, \cite{KJM}, we have discussed spherically symmetric solutions of the Einstein-Euler equations. In this article we study stationary axisymmetric solutions modeling gaseous stars under slow uniform rotation. \\

In the non-relativistic Newtonian theory, the interior structures of gaseous stars are governed by the Euler-Poisson equations. In recent works \cite{JJTM} and \cite{JJTM.Ext} we have constructed uniformly and slowly rotating axisymmetric solutions to the Euler-Poisson equations. In this article we try to extend this result to the relativistic problem. Actually it seems that no existence theorem of the axisymmetric stationary metric in the interior of gaseous
stars in the presence of the pressure, say, neither of vacuum nor of dust, has been established, except only one found in the work by U. Heilig \cite{Heilig}, although many interesting results have been obtained for the vacuum region,
e.g., the Kerr's metric and so on. See e.g. \cite{Kerr}, \cite{KerrGeo}.\\

Let us recall the result of our works \cite{JJTM} and \cite{JJTM.Ext}
on the non-relativistic problem governed by the Euler-Poisson equations.

Assuming $P=\mathsf{A}\rho^{\gamma}$ with $6/5<\gamma< 2$, we constructed a solution with density distribution
$$\rho=\rho_{\mathsf{O}}\Big(\Theta\Big(\frac{r}{\mathsf{a}}, \zeta; \frac{1}{\gamma-1},\mathsf{b}\Big)\vee 0\Big)^{\frac{1}{\gamma-1}},$$
where $\Theta\vee 0=\max(\Theta, 0)$, with the velocity field
$$\vec{v}=(-\Omega y, \Omega x, 0)^T$$
of the Euler-Poisson equations
\begin{align*}
&\frac{\partial\rho}{\partial t}+(\nabla|\rho\vec{v})=0, \\
&\rho\Big(\frac{\partial\vec{v}}{\partial t}+(\vec{v}|\nabla)\vec{v}\Big)+
\nabla P =-\rho\nabla \Phi, \\
&\triangle\Phi=4\pi\mathsf{G}\rho,
\end{align*}
where $(x,y,z)\in\mathbb{R}^3$, $r=\sqrt{x^2+y^2+z^2}$,
$\zeta=z/r$, and $\Omega$ is a constant, the angular velocity of the uniform rotation. 
The last Poisson equation is replaced by the Newton potential
$$\Phi(\vec{x})=-\mathsf{G}\int
\frac{\rho(\vec{x}')}{|\vec{x}-\vec{x}'|}d\vec{x}'
$$
provided that the density $\rho$ is compactly supported.

Here the function $\Theta(r,\zeta; \frac{1}{\gamma-1}, \mathsf{b})$ is the `distorted Lane-Emden
function' with the following properties:

1) The function $\Theta^{\flat}:(x,y,z)\mapsto
\Theta(r,\zeta;\frac{1}{\gamma-1},\mathsf{b})$ is an equatorially and axially symmetric  $C^2$-function on the domain
$\bar{B}(\Xi_0):=\{ (x,y,z) \in \mathbb{R}^3 |  r\leq \Xi_0\}$; Here $\Xi_0$ is a large positive number, which will be fixed in this article, 
such that $\Xi_0 \geq 2\xi_1(\frac{1}{\gamma-1})$, $\xi_1(\frac{1}{\gamma-1})$ being the zero of the Lane-Emden function $\theta(r;\frac{1}{\gamma-1})$ of index
 $\frac{1}{\gamma-1}$, that is, the solution of
$$-\frac{1}{r^2}\frac{d}{dr}r^2\frac{d\theta}{dr}=(\theta\vee 0)^{\frac{1}{\gamma-1}},
\quad \theta|_{r=0}=1; $$

2) $\Theta(0,\zeta; \frac{1}{\gamma-1},\mathsf{b})=1$ and there is a curve
 $\zeta\in [-1,1]\mapsto
r=\Xi_1(\zeta;\frac{1}{\gamma-1},\mathsf{b})$ such that $\Xi_1(\zeta;\frac{1}{\gamma-1},\mathsf{b}) <2\xi_1(\frac{1}{\gamma-1})$ and
$$0\leq r \leq \Xi_0, 0 < \Theta\Big(r,\zeta;\frac{1}{\gamma-1},\mathsf{b}\Big) \qquad \Leftrightarrow \qquad 
0\leq r<\Xi_1\Big(\zeta;\frac{1}{\gamma-1}, \mathsf{b}\Big).$$
 The positive number $\rho_{\mathsf{O}}$ is the central density, an arbitrary positive constant, and 
the parameters $\mathsf{a},\mathsf{b}$ are defined as 
$$\mathsf{a}=\sqrt{\frac{\mathsf{A}\gamma}{4\pi\mathsf{G}(\gamma-1)}}
\rho_{\mathsf{O}}^{-\frac{2-\gamma}{2}},
\quad
\mathsf{b}=\frac{\Omega^2}{2\pi\mathsf{G}\rho_{\mathsf{O}}}.
$$
We require that $\mathsf{b}$ is sufficiently small, say, $\mathsf{b} \leq \varepsilon_0$.

We note that $u=\Theta^{\flat}$ is the equatorially and axially symmetric solution of the integral equation
\begin{equation}
u=\frac{\mathsf{b}}{4}(x^2+y^2)+\mathcal{G}(u), \label{IE}
\end{equation}
where
\begin{align*}
&\mathcal{G}(u)=\mathcal{K}^{(3)}(u\vee 0)^{\frac{1}{\gamma-1}}
-\mathcal{K}^{(3)}(u\vee 0)^{\frac{1}{\gamma-1}}(O)+1, \\
&\mathcal{K}^{(3)}g(\vec{x})=\frac{1}{4\pi}\int\frac{g(\vec{x}')}{|\vec{x}-\vec{x}'|}d\vec{x}'.
\end{align*}
See \cite[Theorem 1]{JJTM}
and \cite[Theorem 2]{JJTM.Ext}. The existence of $\Theta$ is established by the fact that
the Fr\'{e}chet derivative
$$D\mathcal{G}(u)h=\mathcal{K}^{(3)}\Big[\frac{1}{\gamma-1}(u\vee 0)^{\frac{1}{\gamma-1}-1}h\Big]
-\mathcal{K}^{(3)}\Big[\frac{1}{\gamma-1}(u\vee 0)^{\frac{1}{\gamma-1}-1}h\Big](O)$$
of the operator $\mathcal{G}$ in the Banach space of equatorially and axially symmetric continuous functions on $\bar{B}(\Xi_0)$ enjoys the condition
\begin{equation}
\mbox{Ker}(I-D\mathcal{G}(u))=\{0\} \label{HL}
\end{equation}
at $u=\theta(r;\frac{1}{\gamma-1})$, the non-rotating spherically symmetric Lane-Emden 
function. Based on this fact the implicit function theorem guarantees, near $\theta$, the existence of the solution of \eqref {IE}, $\Theta(r,\zeta;\frac{1}{\gamma-1},\mathsf{b})$, for which the condition 
\eqref{HL} holds, too. This property will play an important r\^{o}le in \S 4.3.
\begin{Remark}\label{Remark1}
Here we can take $\Xi_0 >2\xi_1(\frac{1}{\gamma-1})$ arbitrarily large, 
 but we should take $\varepsilon_0$ sufficiently small for the large $\Xi_0$. 
In fact $\Theta$ is the solution of
the equation  \eqref{IE}. We see that the second term of the right-hand side is 
$\mbox{Const.}+O(\frac{1}{r})$ as $r \rightarrow +\infty$ uniformly with respect to $\mathsf{b}$. Therefore we should take $\mathsf{b}\leq\varepsilon_0$ small for large $\Xi_0$ when we want to have $\Theta <0$ for $(\Xi_1(\zeta)<)2\xi_1\leq r<\Xi_0$, by taking into account the growth of the first term $\frac{\mathsf{b}}{4}(x^2+y^2)$. Otherwise, we may cut off the term
$\frac{\mathsf{b}}{4}(x^2+y^2)$ on $ r>2\xi_1$.
\end{Remark}

Anyway let us fix such a solution and denote $\rho=\rho_{\mathsf{N}}$ with $\rho_{\mathsf{O}}=\rho_{\mathsf{N}\mathsf{O}} $ and
put
\begin{equation}
u_{\mathsf{N}}=u_{\mathsf{O}}\Theta\Big(\frac{r}{\mathsf{a}},\zeta;\frac{1}{\gamma-1},
\mathsf{b}\Big)\quad\mbox{with}\quad
u_{\mathsf{O}}=\frac{\mathsf{A}\gamma}{\gamma-1}\rho_{\mathsf{N}\mathsf{O}}^{\gamma-1}.
\end{equation}
Let us denote by $\Phi_{\mathsf{N}}$ the gravitational potential of the density distribution $\rho_{\mathsf{N}}$, that is,
\begin{align*}
\Phi_{\mathsf{N}}(\varpi, z)=&-\mathsf{G}\int_{-\infty}^{+\infty}
\int_0^{+\infty}\int_{0}^{2\pi}
\frac{d\phi'}{\sqrt{\varpi^2+\varpi'^2-2\varpi\varpi'\cos\phi'+(z-z')^2}}\\
&\times \rho_{\mathsf{N}}(\varpi',z')\varpi' d\varpi' dz',
\end{align*}
where $$\varpi=\sqrt{x^2+y^2}, \quad x=\varpi\cos\phi, \quad y=\varpi\sin\phi$$. \\

We are looking for axisymmetric solutions  of the relativistic 
Einstein-Euler equations which approach to this solution of the Euler-Poisson equations as $\mathsf{c}\rightarrow +\infty$.\\

We can describe the main conclusion of this study as follows:

{\bf We are looking for a metric of the form
\begin{equation}
ds^2=e^{2F'}(\mathsf{c}dt+A'd\phi')^2
-e^{-2F'}[e^{2K'}(d\varpi^2+dz^2)+\Pi^2(d\phi')^2],
\end{equation}
where $F', A', K', \Pi$ and the density distribution $\rho$ depend only on $\varpi$ and $z$. Here the co-ordinates are 
$x^0=\mathsf{c}t, x^1=\varpi, x^2=\phi, x^3=z$, and $\phi'=\phi-\Omega t$, while $\Omega$ is the constant, which determined the solution of the 
non-relativistic equations above fixed. The 4-velocity field looked for is
$$U^{\mu}\frac{\partial}{\partial x^{\mu}}=\frac{1}{\mathsf{c}e^{F'}}\Big(\frac{\partial}{\partial t}+
\Omega\frac{\partial}{\partial\phi}\Big).$$
Let us fix an arbitrarily large positive number $R$ with $2\mathsf{a} \xi_1<R\leq\mathsf{a} \Xi_0$. 
If  $u_{\mathsf{O}}/\mathsf{c}^2$ is sufficiently small, then we can construct
 the metric and the density distribution 
on the domain $\mathfrak{D}=\{ (\varpi,z)\  |\  r=\sqrt{\varpi^2+z^2}<R \}$ such that
\begin{align*}
&F'=
-\frac{\Omega^2}{2\mathsf{c}^2}\varpi^2+\frac{\Phi_{\mathsf{N}}}{\mathsf{c}^2}+O\Big(\frac{u_{\mathsf{O}}^2}{\mathsf{c}^4}\Big),
\quad
K'=-\frac{\Omega^2}{2\mathsf{c}^2}{\varpi^2}+O\Big(
\frac{u_{\mathsf{O}}^2}{\mathsf{c}^4}\Big), \\
&A'=-\frac{\Omega}{\mathsf{c}}\varpi^2\Big(1+O\Big(
\frac{u_{\mathsf{O}}}{\mathsf{c}^2}\Big)\Big),
\quad \Pi=\varpi\Big(1+O\Big(\frac{u_{\mathsf{O}}^2}{\mathsf{c}^4}\Big)\Big),
\end{align*}
and
$$\rho=\Big(\frac{\gamma-1}{\mathsf{A}\gamma}\Big)^{\frac{1}{\gamma-1}}
(u\vee 0)^{\frac{1}{\gamma-1}}\Big(1+O\Big(\frac{u}{\mathsf{c}^2}\Big)\Big)\quad\mbox{with}\quad
u=u_{\mathsf{N}}+O\Big(\frac{u_{\mathsf{O}}^2}{\mathsf{c}^2}\Big).
$$
Of course we take $u$ so that $u(0,0)=u_{\mathsf{N}}(0,0)=u_{\mathsf{O}}$. Actually the interior of the star $\{\rho>0\}$ is of the form
$\{(\varpi, z) | r < R(\zeta)\}$
determined by a suitable curve $r=R(\zeta)=\mathsf{a}(\Xi_1(\zeta)+O(u_{\mathsf{O}}^2/\mathsf{c}^2))$.
( Of course $r=\sqrt{\varpi^2+z^2}, \zeta=z/r$.)}\\

The plan of this article is as follows. In \S 2 we will describe how to derive the set of equations to be solved from the Einstein's equations.
Although the set of equations are the same as those found in the book by R. Meinel et al, \cite[(1.34), (1.35)]{Meinel}, \cite{Meinel} gives neither 
how to derive them nor under which conditions they can restore the Einstein's equations. The conditions {\bf (B1), (B2)} are proposed as conditions to restore the original full set of Einstein's equations. Moreover we will discuss on a delicate situation, which appears in quadrature of a first order system for $K'$ when we deal with non-zero pressure field. In usual case one deals with the vacuum or the dust with zero pressure, so Lemma 1 may be a new observation obtained through a tedious calculation. In \S 3 we will derive the set of reduced variables and equations, which are introduced through a heuristic observation of the problem in view of the so called `post-Newtonian approximation'. Little bit of complicated and troublesome calculations are necessary to make clear the estimates needed for the proof of the existence of solutions. \S 4 is devoted to preparation of notations and definitions of functional spaces to be used for it. \S 5 is devoted to the proof of existence of solutions. The problem is divided to two steps. In the first step, the partial differential equations 
for unknown $F', A', \Pi$ to be solved provided that $K'$ is given temporarily are transformed to integral equations involving Newtonian potentials of dimensions 3,4 and 5. In the second step Lemma 2, which is a paraphrase of Lemma 1, will be mobilized to perform the quadrature of the first order system for unknown $K'$. In both steps the fixed point theorem with contraction in suitable Banach spaces prepared in \S 4 can be applied.  In \S 6 we shall clarify the property of the vacuum boundary of the constructed solution with
a compact density support. Finally in \S 7 there will be proposed an open problem on which a continued study should be made. It is the so-called `matter-vacuum matching problem'. This is very important but is still open. 


\section{Basic equations}
This section is devoted to derivation of the basic equations to be solved.\\

We take the coordinates
\begin{equation}
x^0=\mathsf{c}t,\quad x^1=\varpi,\quad x^2=\phi,\quad x^3=z.
\end{equation}

Let us write the metric $ds^2=g_{\mu\nu}dx^{\mu}dx^{\nu}$ in the following form (Lewis 1932 \cite{Lewis}, 
Papapetrou 1966 \cite{Papapetrou} ):
\begin{equation}
ds^2=e^{2F}
(\mathsf{c}dt+Ad\phi)^2-
e^{-2(F-K)}
(d\varpi^2+dz^2)-e^{-2F}\Pi^2d\phi^2, \label{2}
\end{equation}
where the functions $F, A, K$ and $\Pi$ depend only on $\varpi$ and $z$. 
The expression \eqref{2} is called `Lanczos form' after \cite{Lanczos}. 



The components of the metric are as following ( other $g_{\mu\nu}, g^{\mu\nu}$ are zero ):

\begin{align*}
&g_{00}=e^{2F}, \quad
g_{02}=g_{20}=e^{2F}A, \\
&g_{11}=g_{33}=-e^{-2F+2K}, \quad
g_{22}=e^{2F}A^2-e^{-2F}\Pi^2;
\end{align*}
\begin{align*}
& g^{00}=-\frac{1}{\Pi^2}(e^{2F}A^2-e^{-2F}\Pi^2), \quad
 g^{02}=g^{20}=\frac{e^{2F}}{\Pi^2}A, \\
& g^{11}=g^{33}=-e^{2F-2K}, \quad
 g^{22}=-\frac{e^{2F}}{\Pi^2}.
\end{align*}

The stationary and rigid rotation of the fluid is characterized by the 4-velocity field $U^{\mu}$ such that
\begin{equation}
U^0=e^{-G},\quad U^1=U^3=0,\quad
U^2=e^{-G}\frac{\Omega}{\mathsf{c}}.
\end{equation}
Since $U^{\mu}U_{\mu}=1$, the factor $e^{-G}=U^0$ is given by
\begin{equation}
e^{2F}\Big(1+\frac{\Omega}{\mathsf{c}}A\Big)^2
-e^{-2F}\frac{\Omega^2}{\mathsf{c}^2}\Pi^2=e^{2G}. \label{18}
\end{equation}

We shall use the `corotating co-ordinate system' characterized by
\begin{equation}
t'=t,\quad \varpi'=\varpi,\quad
\phi'=\phi-\Omega t, \quad  z'=z, 
\end{equation}
by which the 4-velocity field $U^{\mu'}$ turns out to be
\begin{equation}
U^{0'}=e^{-G},\quad U^{1'}=U^{2'}=U^{3'}=0.
\end{equation}
The metric \eqref{2} can be written of the same form
\begin{equation}
ds^2=e^{2F'}(\mathsf{c}dt'+A'd\phi')^2
-e^{-2(F'-K')}((d\phi')^2+(dz')^2)-e^{-2F'}(\Pi')^2(d\phi')^2 \label{B1}
\end{equation}
with
\begin{equation}
F'-K'=F-K,\quad \Pi'=\Pi, \label{BK}
\end{equation}
provided that
\begin{equation}
e^{2F}\Big(1+\frac{\Omega}{\mathsf{c}}A\Big)^2-
e^{-2F}\frac{\Omega^2}{\mathsf{c}^2}\Pi^2 >0. \label{XP}
\end{equation}
Actually we put
\begin{align}
&e^{2F'}=e^{2F}\Big(1+\frac{\Omega}{\mathsf{c}}A\Big)^2-
e^{-2F}\frac{\Omega^2}{\mathsf{c}^2}\Pi^2, \label{B2} \\
&e^{2F'}A'=e^{2F}\Big(1+\frac{\Omega}{\mathsf{c}}A\Big)A-
e^{-2F}\Pi^2, \label{B3}
\end{align}
provided that \eqref{XP} holds, that is, the right hand side of \eqref{B2} is positive.
Then it can be verified that
\begin{equation}
e^{2F'}(A')^2-e^{-2F'}\Pi^2=
e^{2F}A^2-e^{-2F}\Pi^2
\end{equation}
holds, and therefore \eqref{B1} holds. 
Moreover it can be verified 
\begin{equation}
\Big(1-\frac{\Omega}{\mathsf{c}}A'\Big)e^{2F'}=\Big(1+
\frac{\Omega}{\mathsf{c}}A\Big)e^{2F}. \label{X8}
\end{equation}

Note that comparing \eqref{18} and
\eqref{B2}, we have
\begin{equation}
G=F'.
\end{equation}
Therefore, if we try to find $F', K', A', \Pi$ as functions of $\varpi, z$, then it should be satisfied \\

{\bf (B1)}:\ $\qquad\qquad \displaystyle  
e^{2F'}\Big(1-\frac{\Omega}{\mathsf{c}}A'\Big)^2
-e^{-2F'}\frac{\Omega^2}{\mathsf{c}^2}\Pi^2 >0.$\\

\noindent Only if so, we can define $F, K, A$ by
\begin{align}
&e^{2F}=e^{2F'}\Big(1-\frac{\Omega}{\mathsf{c}}A'\Big)^2
-e^{-2F'}\frac{\Omega^2}{\mathsf{c}^2}\Pi^2, \\
&e^{2F}A=e^{2F'}
\Big(1-\frac{\Omega}{\mathsf{c}}A'\Big)A'
-e^{-2F}\frac{\Omega^2}{\mathsf{c}^2}\Pi^2,
\end{align}
to recover \eqref{2}.\\

Now the energy-momentum tensor of a perfect fluid is
\begin{equation}
T^{\mu\nu}=(\mathsf{c}^2\rho +P)U^{\mu}U^{\nu}-Pg^{\mu\nu}.\label{T}
\end{equation}

The non-zero components of the tensor $T^{\mu\nu}$ are:
\begin{align*}
&T^{00}=e^{-2G}(\mathsf{c}^2\rho+P)+
\frac{P}{\Pi^2}(e^{2F}A^2-e^{-2F}\Pi^2), \\
&T^{02}=T^{20}=
e^{-2G}(\mathsf{c}^2\rho+P)\frac{\Omega}{\mathsf{c}}-
\frac{P}{\Pi^2}e^{2F}A, \\
&T^{11}=T^{33}=Pe^{2F-2K}, \\
&T^{22}=
e^{-2G}(\mathsf{c}^2\rho+P)\frac{\Omega^2}{\mathsf{c}^2}+
\frac{P}{\Pi^2}e^{2F}.
\end{align*}

Keeping in mind \eqref{18}, we can show through tedious calculations that 
$$\nabla_{\mu}T^{\mu 1}=
e^{2F-2K}(\partial_1P+(\mathsf{c}^2\rho+P)\partial_1G) $$
and
$$\nabla_{\mu}T^{\mu 3}=
e^{2F-2K}(\partial_3P+(\mathsf{c}^2\rho+P)\partial_3G). $$
Hence the identity $\nabla_{\mu}T^{\mu 1}=\nabla_{\mu}T^{\mu 3}=0$ reduces
\begin{equation}
\partial_1P+(\mathsf{c}^2\rho+P)\partial_1G=0,\qquad
\partial_3P+(\mathsf{c}^2\rho+P)\partial_3G=0.
\end{equation}

Here and hereafter $\partial_{\mu}$ stands for $\displaystyle \frac{\partial}{\partial x^{\mu}}$, and 
$\nabla_{\lambda}$ denotes the covariant derivative:
$$\nabla_{\lambda}A^{\mu\nu}=
\partial_{\lambda}A^{\mu\nu}+
\Gamma_{\alpha\lambda}^{\mu}A^{\alpha\nu}+
\Gamma_{\alpha\lambda}^{\nu}A^{\mu\alpha}, $$
where $\Gamma_{\nu\lambda}^{\mu}$ is the Christoffel symbol:
$$\Gamma_{\nu\lambda}^{\mu}=
\frac{1}{2}g^{\mu\alpha}(\partial_{\lambda}g_{\alpha\nu}+\partial_{\nu}g_{\alpha\lambda}
-\partial_{\alpha}g_{\nu\lambda}). $$

Therefore, defining the `relativistic enthalpy density ' $u$ by
\begin{equation}
u:=\mathsf{c}^2\int_0^{\rho}\frac{dP}{\mathsf{c}^2\rho+P}=
\int_0^{\rho}\frac{dP}{\rho+P/\mathsf{c}^2},
\end{equation}
we have
\begin{equation}
\frac{u}{\mathsf{c}^2}+G=\mbox{Const.},
\end{equation}
or
\begin{equation}
F'=G=-\frac{u}{\mathsf{c}^2}+\mbox{Const.}
\end{equation}

Note that the other Euler equations $$\nabla_{\mu}T^{\mu 0}=0, \quad
\nabla_{\mu}T^{\mu 2}=0 $$
hold automatically. \\

We are going to derive the Einstein's equations for the metric \eqref{2} and the energy-momentum tensor $T^{\mu\nu}$ \eqref{T} in the form \cite[(1. 34)(1.35)]{Meinel} by following the calculations described in 
\cite{Islam}.\\

In order to write down the Einstein's equations
\begin{equation}
R_{\mu\nu}-\frac{1}{2}g_{\mu\nu}R=
\frac{8\pi\mathsf{G}}{\mathsf{c}^4}T_{\mu\nu}
\end{equation}
or
\begin{equation}
R_{\mu\nu}=\frac{8\pi\mathsf{G}}{\mathsf{c}^4}
(T_{\mu\nu}-\frac{1}{2}g_{\mu\nu}T),
\end{equation}
where $R$ stands for $g^{\alpha\beta}R_{\alpha\beta}$ and
$T$ stands for $g_{\alpha\beta}T^{\alpha\beta}$, let us compute $U_{\mu}$,
$T_{\mu\nu}$ and $T$. The result is as following:

\begin{align*}
&U_0=e^{2F-G}\Big(1+A\frac{\Omega}{\mathsf{c}}\Big), \\
&U_1=U_3=0, \\
&U_2=e^{-G}\Big(
e^{2F}A+(e^{2F}A^2-e^{-2F}\Pi^2)\frac{\Omega}{\mathsf{c}}\Big) ;
\end{align*}
\begin{align*}
T_{00}&=(\mathsf{c}^2\rho+P)e^{4F-2G}\Big(1+A\frac{\Omega}{\mathsf{c}}\Big)^2
-Pe^{2F}, \\
T_{02}&=T_{20}= \\
&=(\mathsf{c}^2\rho+P)e^{2F-2G}\Big(1+A\frac{\Omega}{\mathsf{c}}\Big)
\Big(e^{2F}A+(e^{2F}A^2-
e^{-2F}\Pi^2)\frac{\Omega}{\mathsf{c}}\Big)
-Pe^{2F}A, \\
T_{11}&=T_{33}=Pe^{-2F+2K}, \\
T_{22}&=
(\mathsf{c}^2\rho+P)e^{-2G}
\Big(e^{2F}A+(e^{2F}A^2-e^{-2F}\Pi^2)\frac{\Omega}{\mathsf{c}}\Big)^2
-P(e^{2F}A^2-e^{-2F}\Pi^2)
\end{align*}
and other $T_{\mu\nu}$'s are zero; We have
\begin{align*}
T&=(\mathsf{c}^2\rho+P)e^{-2G}\Big[e^{2F}\Big(1+A\frac{\Omega}{\mathsf{c}}\Big)^2
-e^{-2F}\Pi^2\frac{\Omega^2}{\mathsf{c}^2}\Big]
-4P \\
&=\mathsf{c}^2\rho-3P.
\end{align*}\\


In order to calculate the Christoffel symbols
$\Gamma_{\nu\lambda}^{\mu}$
and the Ricci tensor
$$R_{\mu\nu}=\partial_{\alpha}\Gamma_{\mu\nu}^{\alpha}
-\partial_{\nu}\Gamma_{\mu\alpha}^{\alpha}
+\Gamma_{\mu\nu}^{\alpha}\Gamma_{\alpha\beta}^{\beta}
-\Gamma_{\mu\alpha}^{\beta}\Gamma_{\nu\beta}^{\alpha} $$ 
 and to write down explicitly the Einstein's equations, it is convenient to write the metric as
\begin{equation}
ds^2=f\mathsf{c}^2dt^2-2k\mathsf{c}dtd\phi-ld\phi^2
-e^m(d\varpi^2+dz^2), \label{25}
\end{equation}
that is, to put
\begin{equation}
f=e^{2F},\quad k=-e^{2F}A, \quad m=2(-F+K),\quad l=-e^{2F}A^2+
e^{-2F}\Pi^2.
\end{equation}
The expression \eqref{25} is called `Lewis metric' after \cite{Lewis}\\

First let us note the identity
\begin{equation}
\Pi^2=fl+k^2. \label{27}
\end{equation}
and that the factor $e^{-G}=U^0$ is given by
\begin{equation}
(U^0)^{-2}=e^{2G}=f-2\frac{\Omega}{\mathsf{c}}k
-\frac{\Omega^2}{\mathsf{c}^2}l.
\end{equation}

Let us introduce the quantity $\Sigma$ defined by
\begin{equation}
\Sigma=\partial_1f\cdot\partial_1l+
\partial_3f\cdot\partial_3l+(\partial_1k)^2+(\partial_3k)^2.
\end{equation}\\

Now the components of the metric and the Christoffel symbols are as following (other
$g_{\mu\nu}, g^{\mu\nu}, \Gamma_{\mu\nu}^{\lambda}$ are zero ):
$$
g_{00}=f, \quad
g_{02}=g_{20}=-k, \quad g_{22}=-l, \quad g_{11}=g_{33}=-e^m;
$$
$$
g^{00}=\frac{l}{\Pi^2}, \quad g^{02}=g^{20}=-\frac{k}{\Pi^2}, \quad g^{22}=-\frac{f}{\Pi^2}, \quad g^{11}=g^{33}=-e^{-m};
$$
\begin{align*}
&\Gamma_{01}^0=\Gamma_{10}^0=\frac{1}{2\Pi^2}(l\partial_1f+k\partial_1k), \quad \Gamma_{03}^0=\Gamma_{30}^0=\frac{1}{2\Pi^2}(l\partial_3f+k\partial_3k), \\
&\Gamma_{12}^0=\Gamma_{21}^0=\frac{1}{\Pi^2}(k\partial_1l-l\partial_1k), \quad \Gamma_{23}^0=\Gamma_{32}^0=\frac{1}{2\Pi^2}(k\partial_3l-l\partial_3k);
\end{align*}
\begin{align*}
&\Gamma_{00}^1=\frac{1}{2}e^{-m}\partial_1f, \quad \Gamma_{02}^1=\Gamma_{20}^1=-\frac{1}{2}e^{-m}\partial_1k, \quad \Gamma_{11}^1=\frac{1}{2}\partial_1m, \\
&\Gamma_{13}^1=\Gamma_{31}^1=\frac{1}{2}\partial_3m, \quad \Gamma_{33}^1=-\frac{1}{2}\partial_1m, \quad \Gamma_{22}^1=-\frac{1}{2}e^{-m}\partial_1l;
\end{align*}
\begin{align*}
&\Gamma_{01}^2=\Gamma_{10}^2=\frac{1}{2\Pi^2}
(f\partial_1k-k\partial_1f), \quad\Gamma_{03}^2=\Gamma_{30}^2=\frac{1}{2\Pi^2}
(f\partial_3k-k\partial_3f), \\
&\Gamma_{12}^2=\Gamma_{21}^2=\frac{1}{2\Pi^2}
(f\partial_1l+k\partial_1k), \quad\Gamma_{23}^2=\Gamma_{32}^2=\frac{1}{2\Pi^2}
(f\partial_3l+k\partial_3k);
\end{align*}
\begin{align*}
&\Gamma_{00}^3=\frac{1}{2}e^{-m}\partial_3f, \quad\Gamma_{02}^3=\Gamma_{20}^3=-\frac{1}{2}e^{-m}\partial_3k, \quad\Gamma_{11}^3=-\frac{1}{2}\partial_3m, \\
&\Gamma_{13}^3=\Gamma_{31}^3=\frac{1}{2}\partial_1m, \quad\Gamma_{33}^3=\frac{1}{2}\partial_3m, \quad\Gamma_{22}^3=-\frac{1}{2}e^{-m}\partial_3l.
\end{align*}

The components of the Ricci tensor are as following (other $R_{\mu\nu}$'s are zero ):
\begin{subequations}
\begin{align}
\frac{2e^m}{\Pi}R_{00}&=\partial_1\frac{\partial_1f}{\Pi}+\partial_3\frac{\partial_3f}{\Pi}+
\frac{1}{\Pi^3}f\Sigma, \label{30a}\\
-\frac{2e^m}{\Pi}R_{02}&=-\frac{2e^m}{\Pi}R_{20}=
\partial_1\frac{\partial_1k}{\Pi}+\partial_3\frac{\partial_3k}{\Pi}+\frac{1}{\Pi^3}k\Sigma, \label{30b} \\
-\frac{2e^m}{\Pi}R_{22}&=
\partial_1\frac{\partial_1l}{\Pi}+\partial_3\frac{\partial_3l}{\Pi}+\frac{1}{\Pi^3}l\Sigma, \label{30c}\\
2R_{11}&=-\partial_1^2m-\partial_3^2m-2\frac{\partial_1^2\Pi}{\Pi}
+\frac{1}{\Pi}(\partial_1m\partial_1\Pi-\partial_3m\partial_3\Pi) + \nonumber \\
&+\frac{1}{\Pi^2}(\partial_1f\partial_1l+(\partial_1k)^2), \label{30d}\\
2R_{33}&=-\partial_1^2m-\partial_3^2m-2\frac{\partial_3^2\Pi}{\Pi}
-\frac{1}{\Pi}(\partial_1m\partial_1\Pi-\partial_3m\partial_3\Pi) + \nonumber \\
&+\frac{1}{\Pi^2}(\partial_3f\partial_3l+(\partial_3k)^2), \label{30e}\\
2R_{13}&=2R_{31}=-2\frac{\partial_1\partial_3\Pi}{\Pi}+
\frac{1}{\Pi}(\partial_3m\partial_1\Pi+\partial_1m\partial_3\Pi)
\nonumber \\
&+\frac{1}{2\Pi^2}(\partial_1f\partial_3l+\partial_1l\partial_3f
+2\partial_1k\partial_3k).\label{30f}
\end{align}
\end{subequations}
During the calculations we have used the identities
\begin{equation}
\Gamma_{j\mu}^{\mu}=\partial_jm+\frac{\partial_j\Pi}{\Pi} \qquad
\mbox{for}\quad j=1,3.
\end{equation}\\

Now the components of the 4-velocity vector are:
\begin{align*}
&U^0=e^{-G}=\Big(f-
2\frac{\Omega}{\mathsf{c}}k-\frac{\Omega^2}{\mathsf{c}^2}l\Big)^{-1/2},
\quad U^1=U^3=0,\quad U^2=\frac{\Omega}{\mathsf{c}}U^0; \\
&U_0=\Big(f-\frac{\Omega}{\mathsf{c}}k\Big)U^0,\quad
U^1=U^3=0,\quad
U^2=-\Big(k+\frac{\Omega}{\mathsf{c}}\Big)U^0.
\end{align*}

The Einstein's equations are
\begin{equation}
R_{\mu\nu}=\frac{8\pi\mathsf{G}}{\mathsf{c}^4}\mathfrak{T}_{\mu\nu}, \quad\mbox{where}\quad
\mathfrak{T}_{\mu\nu}:=T_{\mu\nu}-\frac{1}{2}g_{\mu\nu}T.
\end{equation}

The components $\mathfrak{T}_{\mu\nu}$ turn out to be as following (other
$\mathfrak{T}_{\mu\nu}$'s are zero) :
\begin{subequations}
\begin{align}
&\mathfrak{T}_{00}=
\frac{1}{2}(\mathsf{c}^2\rho+P)e^{-2G}
\Big[\Big(f-\frac{\Omega}{\mathsf{c}}k\Big)^2+\frac{\Omega^2}{\mathsf{c}^2}\Pi^2\Big]
+Pf, \label{Ta}\\
&\mathfrak{T}_{02}=\mathfrak{T}_{20}=
\frac{1}{2}(\mathsf{c}^2\rho+P)e^{-2G}
\Big[-kf-2\frac{\Omega}{\mathsf{c}}fl+\frac{\Omega^2}{\mathsf{c}^2}kl\Big]-Pk, \label{Tb}\\
&\mathfrak{T}_{22}=
\frac{1}{2}(\mathsf{c}^2\rho+P)e^{-2G}\Big[\Pi^2+
\Big(k+\frac{\Omega}{\mathsf{c}}l\Big)^2\Big]-Pl, \label{Tc}\\
&\mathfrak{T}_{11}=\mathfrak{T}_{33}=\frac{e^m}{2}(\mathsf{c}^2\rho -P).
\end{align}
\end{subequations}

Note that
\begin{equation}
T=\mathsf{c}^2\rho-3P.
\end{equation}\\

Thus the Einstein's equations are
\begin{subequations}
\begin{align}
&R_{00}=\frac{8\pi\mathsf{G}}{\mathsf{c}^4}\mathfrak{T}_{00}, \label{36a}\\
&R_{02}=\frac{8\pi\mathsf{G}}{\mathsf{c}^4}\mathfrak{T}_{02}, \label{36b}\\
&R_{22}=\frac{8\pi\mathsf{G}}{\mathsf{c}^4}\mathfrak{T}_{22}, \label{36c}\\
&R_{11}=\frac{8\pi\mathsf{G}}{\mathsf{c}^4}\mathfrak{T}_{11}, \label{W30d}\\
&R_{33}=\frac{8\pi\mathsf{G}}{\mathsf{c}^4}\mathfrak{T}_{33}, \label{W30e}\\
&R_{13}=R_{31}=0
\end{align}
\end{subequations}\\

Now we put
\begin{equation}
f':=e^{2F'}=e^{2G}=f-
2\frac{\Omega}{\mathsf{c}}k-\frac{\Omega^2}{\mathsf{c}^2}l, \quad
k':=k+\frac{\Omega}{\mathsf{c}}l, \quad
l'=l, \quad m'=m.
\end{equation}
Using \eqref{BK}\eqref{B2}\eqref{X8}, we can verify that 
$$k'=-e^{2F'}A', \qquad m'=2(-F'+K'). $$

Then we claim the following

\begin{Proposition}\label{PropositionA}
The set of equations \eqref{36a}\eqref{36b}\eqref{36c} implies
\begin{subequations}
\begin{align}
&\frac{\partial^2F'}{\partial\varpi^2}+\frac{\partial^2F'}{\partial z^2}+
\frac{1}{\Pi}\Big(\frac{\partial F'}{\partial\varpi}\frac{\partial\Pi}{\partial\varpi}+
\frac{\partial F'}{\partial z}\frac{\partial\Pi}{\partial z}\Big)+
\frac{e^{4F'}}{2\Pi^2}
\Big[\Big(\frac{\partial A'}{\partial\varpi}\Big)^2+\Big(\frac{\partial  A'}{\partial z}\Big)^2\Big] \nonumber \\
&=\frac{4\pi\mathsf{G}}{\mathsf{c}^4}e^{-2F'+2K'}(\mathsf{c}^2\rho+3P),\label{EAa} \\
&\frac{\partial}{\partial\varpi}\Big(\frac{e^{4F'}}{\Pi}\frac{\partial A'}{\partial\varpi}\Big)+
\frac{\partial}{\partial z}\Big(\frac{e^{4F'}}{\Pi}\frac{\partial A'}{\partial z}\Big)=0, \label{EAb} \\
&\frac{\partial^2\Pi}{\partial\varpi^2}+\frac{\partial^2\Pi}{\partial z^2}=\frac{16\pi\mathsf{G}}{\mathsf{c}^4}
e^{-2F'+2K'}P\Pi. \label{EAc}
\end{align}
\end{subequations}
\end{Proposition}

Proof. 
First we note that we have the identity
\begin{equation}
\frac{e^m}{\Pi}(lR_{00}-2kR_{02}-fR_{22})=
\partial_1^2\Pi+\partial_3^2\Pi,\label{32}
\end{equation}
which can be verified from \eqref{30a}\eqref{30b}\eqref{30c} thanks to \eqref{27}.

On the other hand, we have the identity
\begin{equation}
l\mathfrak{T}_{00}-2k\mathfrak{T}_{02}
-f\mathfrak{T}_{22}=2P\Pi^2, \label{37}
\end{equation}
which can be verified from \eqref{Ta}\eqref{Tb}\eqref{Tc} thanks to \eqref{27}. 
Hence \eqref{32} and \eqref{37} leads us to the equation
\eqref{EAc}.

\begin{Remark} If we consider the dust for which $P=0$, then \eqref{EAc} says that $\Pi(\varpi, z)$ is a harmonic function and we can assume $\Pi=\varpi$ by conformal change of coordinates.
See \cite[p. 26]{Islam}. (This was first used by 
\cite{Weyl} and generalized to the present case by
\cite{Lewis}.)  But it is not the case when $P\not=0$. 
\end{Remark}

Next we can verify that
\begin{equation}
\Big(k+\frac{\Omega}{\mathsf{c}}l\Big)\mathfrak{T}_{00}+
\Big(f+\frac{\Omega^2}{\mathsf{c}^2}l\Big)\mathfrak{T}_{02}
+\frac{\Omega}{\mathsf{c}}\Big(f-\frac{\Omega}{\mathsf{c}}k\Big)\mathfrak{T}_{22}=0.\label{W33}
\end{equation}
On the other hand, it can be verified from \eqref{30a}\eqref{30b}\eqref{30c} that
\begin{align}
&-\frac{2e^m}{\Pi}
\Big[\Big(k+\frac{\Omega}{\mathsf{c}}l\Big)R_{00}+
\Big(f+\frac{\Omega^2}{\mathsf{c}^2}l\Big)R_{02}+
\frac{\Omega}{\mathsf{c}}\Big(f-\frac{\Omega}{\mathsf{c}}k\Big)R_{22}\Big] = \nonumber \\
&=\partial_1\frac{1}{\Pi}(f'\partial_1k'-k'\partial_1f')+
\partial_3\frac{1}{\Pi}(f'\partial_3k'-k'\partial_3f').
 \label{W34}
\end{align}

Let us note the identities
\begin{align}
&\Big(k+\frac{\Omega}{\mathsf{c}}l\Big)f
-\Big(f+\frac{\Omega^2}{\mathsf{c}^2}l\Big)k
-\frac{\Omega}{\mathsf{c}}\Big(f-\frac{\Omega}{\mathsf{c}}\Big)l=0, \\
&\Sigma=\partial_1f'\partial_1l'+\partial_3f'\partial_3l'+(\partial_1k')^2
+(\partial_3k')^2, \\
&\Pi^2=f'l'+(k')^2.
\end{align}\\

Anyway 
it follows from \eqref{W33}\eqref{W34} that
\begin{equation}
\partial_1\frac{1}{\Pi}(f'\partial_1k'-k'\partial_1f')+
\partial_3\frac{1}{\Pi}(f'\partial_3k'-k'\partial_3f')=0, \label{45a}
\end{equation}
 Note that
$$f'\partial_jk'-k'\partial_jf'= -e^{4F'}\partial_jA', \qquad (j=1,3).$$
Hence \eqref{45a} reads
\eqref{EAb},
which is nothing but \cite[(1.34b)]{Meinel}.\\

Now it follows from \eqref{Ta}\eqref{Tb}\eqref{Tc} that
\begin{equation}
\mathfrak{T}_{00}+2\frac{\Omega}{\mathsf{c}}\mathfrak{T}_{02}
+\frac{\Omega^2}{\mathsf{c}^2}\mathfrak{T}_{22}=
\frac{1}{2}(\mathsf{c}^2\rho+3P)f'.
\end{equation}
Here we have used the identity
$$e^{2G}=f'=f-\frac{\Omega}{\mathsf{c}}-\frac{\Omega^2}{\mathsf{c}^2}l. $$
On the other hand it follows from \eqref{30a}\eqref{30b}\eqref{30c} that
\begin{equation}
\frac{2e^m}{\Pi}\Big(
R_{00}+2\frac{\Omega}{\mathsf{c}}R_{02}+\frac{\Omega^2}{\mathsf{c}^2}R_{22}\Big)=
\partial_1\frac{\partial_1f'}{\Pi}+\partial_3\frac{\partial_3f'}{\Pi}+\frac{f'}{\Pi^3}\Sigma.
\end{equation}
Therefore the Einstein's equations imply
\begin{equation}
\partial_1\frac{\partial_1f'}{\Pi}+\partial_3\frac{\partial_3f'}{\Pi}+\frac{f'}{\Pi^3}\Sigma=\frac{8\pi\mathsf{G}}{\mathsf{c}^4}
e^m\frac{f'}{\Pi}(\mathsf{c}^2\rho+3P). \label{49}
\end{equation}
But we have the identity
\begin{align}
\Sigma&=2\Pi\partial_1\Pi\frac{\partial_1f'}{f'}-\Big(\frac{\Pi\partial_1f'}{f'}\Big)^2
+(f'\partial_1A')^2 + \nonumber \\
&+
2\Pi\partial_3\Pi\frac{\partial_3f'}{f'}-\Big(\frac{\Pi\partial_3f'}{f'}\Big)^2
+(f'\partial_3A')^2,
\end{align}
which can be verified by using $k'=-f'A'$ and
$$\Sigma=\partial_1f'\partial_1l'+\partial_3f'\partial_3l'+(f'\partial_1A')^2+(f'\partial_3A')^2 $$
with $\Pi^2=f'l'+(k')^2$. Hence \eqref{49} reads
\begin{align}
&\frac{1}{f'}(\partial_1^2f'+\partial_3^2f')+
\frac{1}{f'\Pi}(\partial_1f'\partial_1\Pi+\partial_3f'\partial_3\Pi) \nonumber \\
&-\Big[\Big(\frac{\partial_1f'}{f'}\Big)^2+
\Big(\frac{\partial_3f'}{f'}\Big)^2\Big]
+\Big(\frac{f'}{\Pi}\Big)^2
\Big[(\partial_1A')^2+(\partial_3A')^2\Big]
=\frac{8\pi\mathsf{G}}{\mathsf{c}^4}e^m(\mathsf{c}^2\rho +3P).
\end{align}
This is nothing but \cite[(1.34a)]{Meinel}, that is,
we have \eqref{EAa},
since $2F'=\log f'$ so that
$$2\partial_jF'=\frac{\partial_jf'}{f'}, \qquad
2\partial_j^2F'=\frac{\partial_j^2f'}{f'}-
\Big(\frac{\partial_jf'}{f'}\Big)^2,\quad j=1,3.
$$

Proof of Proposition~\ref{PropositionA} is finished.\vspace{6mm} $\square$\\

Inversely we can show that \eqref{EAa}\eqref{EAb}\eqref{EAc} imply 
\eqref{36a}\eqref{36b}\eqref{36c}.
In fact, if we put
$$Q_{\mu\nu}:=R_{\mu\nu}-\frac{8\pi\mathsf{G}}{\mathsf{c}^4}\mathfrak{T}_{\mu\nu},
$$
then \eqref{36a}\eqref{36b}\eqref{36c} claim that $Q_{00}=Q_{02}=Q_{22}=0$.
On the other hand \eqref{EAa}\eqref{EAb}\eqref{EAc} are nothing but
\begin{align*}
&lQ_{00}-2kQ_{02}-fQ_{22}=0, \\
&\Big(k+\frac{\Omega}{\mathsf{c}}l\Big)Q_{00}
+\Big(f+\frac{\Omega^2}{\mathsf{c}^2}l\Big)Q_{02}+
\frac{\Omega}{\mathsf{c}}\Big(f-\frac{\Omega}{\mathsf{c}}k\Big)Q_{22}=0, \\
&Q_{00}+2\frac{\Omega}{\mathsf{c}}Q_{02}+\frac{\Omega^2}{\mathsf{c}^2}Q_{22}=0.
\end{align*}
Here we see 
\begin{align}
\Delta&:=\det
\begin{bmatrix}
l & -2k & -f \\
k+\frac{\Omega}{\mathsf{c}}l & f+\frac{\Omega^2}{\mathsf{c}^2}l & \frac{\Omega}{\mathsf{c}}\Big(f-\frac{\Omega}{\mathsf{c}}k\Big) \\
1 & 2\frac{\Omega}{\mathsf{c}} & \frac{\Omega^2}{\mathsf{c}^2} 
\end{bmatrix} \nonumber \\
&=\Big(2k+\frac{\Omega}{\mathsf{c}}l\Big)^2\frac{\Omega^2}{\mathsf{c}^2}
+f(f-\frac{4\Omega}{\mathsf{c}}k-\frac{2\Omega^2}{\mathsf{c}^2}l\Big) \nonumber \\
&=(f')^2=e^{4F'}\not= 0.
\end{align}

Thus we can claim

\begin{Proposition}\label{PropositionB}
The set of equations
\eqref{EAa}, \eqref{EAb}, \eqref{EAc} implies 
\eqref{36a}, \eqref{36b}, \eqref{36c}.
\end{Proposition}

Proof. $\Delta \not=0$ guarantees
 $Q_{00}=Q_{02}=Q_{22}=0$.
This means \eqref{36a},\eqref{36b},\eqref{36c}. $\square$. \\

Next we claim

\begin{Proposition}\label{PropositionC}
The equations \eqref{W30d}\eqref{W30e} imply
\begin{subequations}
\begin{align}
&\frac{\partial\Pi}{\partial\varpi}\frac{\partial K'}{\partial\varpi}-
\frac{\partial\Pi}{\partial z}\frac{\partial K'}{\partial z}=
\frac{1}{2}\Big(\frac{\partial^2\Pi}{\partial \varpi^2}-
\frac{\partial^2\Pi}{\partial z^2}\Big)+
\Pi\Big[\Big(\frac{\partial F'}{\partial\varpi}\Big)^2-
\Big(\frac{\partial F'}{\partial z}\Big)^2\Big] + \nonumber \\
&-\frac{e^{4F'}}{4\Pi}
\Big[\Big(\frac{\partial A'}{\partial\varpi}\Big)^2-
\Big(\frac{\partial A'}{\partial z}\Big)^2\Big],\label{ECd}\\
&\frac{\partial\Pi}{\partial z}\frac{\partial K'}{\partial\varpi}
+\frac{\partial\Pi}{\partial\varpi}\frac{\partial K'}{\partial z}=
\frac{\partial^2\Pi}{\partial \varpi \partial z}+2\Pi
\frac{\partial F'}{\partial\varpi}\frac{\partial F'}{\partial z}-
\frac{e^{4F'}}{2\Pi}
\frac{\partial A'}{\partial\varpi}
\frac{\partial A'}{\partial z},\label{ECe}
\end{align}
\end{subequations}
\end{Proposition}

Proof.
It follows from \eqref{30d}\eqref{30e} that
\begin{align*}
2(R_{11}-R_{33})&=-\frac{2}{\Pi}(\partial_1^2\Pi-\partial_3^2\Pi)+
\frac{2}{\Pi}(\partial_1m\partial_1\Pi-\partial_3m\partial_3\Pi)+ \\
&+\frac{1}{\Pi^2}
(\partial_1f\partial_1l+(\partial_1k)^2
-\partial_3f\partial_3l-(\partial_3k)^2).
\end{align*}
On the other hand, since $\mathfrak{T}_{11}=\mathfrak{T}_{33}$, 
the Einstein's equations imply
\begin{align}
-\frac{2}{\Pi}(\partial_1^2\Pi-\partial_3^2\Pi)+&
\frac{2}{\Pi}(\partial_1m\partial_1\Pi-\partial_3m\partial_3\Pi)+ \nonumber \\
&+\frac{1}{\Pi^2}
(\partial_1f\partial_1l+(\partial_1k)^2
-\partial_3f\partial_3l-(\partial_3k)^2)=0.\label{53}
\end{align}
Using $k'=-f'A', \Pi=f'l'+(k')^2$, we can show that
\begin{align}
\partial_jf\partial_jl+(\partial_jk)^2&=
\partial_jf'\partial_jl'+(\partial_jk')^2 \nonumber\\
&=-\Big(\frac{\partial_jf'}{f'}\Big)^2\Pi^2+
2\frac{\partial_jf'}{f'}\Pi\partial_j\Pi+
(f')^2(\partial_jA')^2 \nonumber \\
&=-4(\partial_jF')^2\Pi^2
+4\partial_jF'\Pi\partial_j\Pi+
e^{4F'}(\partial_jA')^2.
\end{align}
Therefore, since $m=m'=-2F'+2K'$, \eqref{53} reads
\eqref{ECd},
which is nothing but \cite[(1.35a)]{Meinel}.\\

Since $\mathfrak{T}_{13}=0$, the Einstein's equations imply that $R_{13}=0$. Hence it follows from \eqref{30f} that
\begin{align*}
-2\frac{\partial_1\partial_3\Pi}{\Pi}&+
\frac{1}{\Pi}(\partial_3m\partial_1\Pi+
\partial_1m\partial_3\Pi) + \\
&+\frac{1}{2\Pi^2}
(\partial_1f\partial_3l+\partial_1l\partial_3f+2\partial_1k\partial_3k)=0.
\end{align*}
But we see that
\begin{align}
&\partial_1f\partial_3l+\partial_1l\partial_3f+2\partial_1k\partial_3k=
\partial_1f'\partial_3l'+\partial_1l'\partial_3f'+2\partial_1k'\partial_3k' \nonumber \\
&=-2\frac{\partial_1f'\partial_3f'}{(f')^2}\Pi^2+
\frac{2\Pi}{f'}(\partial_1f'\partial_3\Pi+\partial_3f'\partial_1\Pi)
+2(f')^2\partial_1A'\partial_3A' \nonumber \\
&=-8\partial_1F'\partial_3F'\Pi^2+
4\Pi(\partial_1F'\partial_3\Pi+\partial_3F'\partial_1\Pi)+
2e^{4F'}\partial_1A'\partial_3A'.
\end{align}
Therefore, since $m=m'=-2F'+2K'$, $R_{13}=0$ reads
\eqref{ECe},
which is nothing but \cite[(1.35b)]{Meinel}. \vspace{6mm} $\square$\\

In order to consider the inverse, we put the following assumption, keeping in mind the limiting case $\Pi=\varpi$:\\

{\bf (B2): \  $\displaystyle \Big(\frac{\partial\Pi}{\partial\varpi}\Big)^2+\Big(\frac{\partial\Pi}{\partial z}\Big)^2$ does not vanish. }\\

Then we can claim the following

\begin{Proposition}\label{PropositionD}
If the assumption
{\bf (B2)} holds, then \eqref{ECd}\eqref{ECe} imply \eqref{W30d}\eqref{W30e},
provided that \eqref{36a}\eqref{36b}\eqref{36c} hold.
\end{Proposition}

Proof. Suppose
that \eqref{ECd},\eqref{ECe} hold. Then we have
$R_{11}=R_{33}$ together with $R_{13}=0$. Let us consider
$$ Q^{\mu\nu}=R^{\mu\nu}-\frac{8\pi\mathsf{G}}{\mathsf{c}^4}\mathfrak{T}^{\mu\nu}.$$
We already know that $Q^{\mu\nu}=0$ for $(\mu,\nu)\not=(1,1),(3,3)$ and $Q^{11}=Q^{33}$. 
We want to show that $\mathcal{Q}:=Q^{11}=Q^{33}$ vanishes.

Since the Euler equations $\nabla_{\mu}T^{\mu\nu}=0$ and the Bianchi identities
 
\noindent $\nabla_{\mu}(R^{\mu\nu}-\frac{1}{2}g^{\mu\nu}R)=0$ hold, we have
$$\nabla_{\mu}Q^{\mu j}=
\frac{1}{2}g^{jj}\partial_j\Big(R+\frac{8\pi\mathsf{G}}{\mathsf{c}^4}T\Big)=
-\frac{1}{2}e^{-m}\partial_j\Big(R+\frac{8\pi\mathsf{G}}{\mathsf{c}^4}T\Big),
\quad j=1,3.$$
But
$$\nabla_{\mu}Q^{\mu j}=\Big[\partial_j+\Big(\partial_jm+\frac{\partial_j\Pi}{\Pi}\Big)\Big]\mathcal{Q}, 
\quad j=1,3,$$
since we already know $Q^{\mu\nu}=0$ for $(\mu,\nu)\not= (1,1), (3,3)$
and 
$$\Gamma_{j\mu}^{\mu}=\partial_jm+\frac{\partial_j\Pi}{\Pi}, \quad\mbox{and}\quad
\Gamma_{11}^j+\Gamma_{33}^j=0 \quad\mbox{for}\quad j=1,3.$$
On the other hand, by contraction we have
$$
(g_{11}+g_{33})\mathcal{Q}=
-2e^m\mathcal{Q}=R+\frac{8\pi\mathsf{G}}{\mathsf{c}^4}T.$$
Therefore
$$\Big[\partial_j+\Big(\partial_jm+\frac{\partial_j\Pi}{\Pi}\Big)\Big]\mathcal{Q}=-\frac{1}{2}e^{-m}\partial_j
(-2e^m\mathcal{Q})=[\partial_j+(\partial_jm)]\mathcal{Q}, 
$$
for $j=1,3$,
that is, $\displaystyle \frac{\partial_j\Pi}{\Pi}\mathcal{Q}=0$ for $j=1,3$.
Under the assumption {\bf (B2)}, it should hold that $\mathcal{Q}=Q^{11}=Q^{33}=0$. So, 
under this assumption, we can claim that \eqref{W30d}, \eqref{W30e} hold. \hspace{8mm} $\square$\\

Summing up, the system of equations to be solved turns out to be
\begin{subequations}
\begin{align}
&\frac{\partial^2F'}{\partial\varpi^2}+\frac{\partial^2F'}{\partial z^2}+
\frac{1}{\Pi}\Big(\frac{\partial F'}{\partial\varpi}\frac{\partial\Pi}{\partial\varpi}+
\frac{\partial F'}{\partial z}\frac{\partial\Pi}{\partial z}\Big)+
\frac{e^{4F'}}{2\Pi^2}
\Big[\Big(\frac{\partial A'}{\partial\varpi}\Big)^2+\Big(\frac{\partial  A'}{\partial z}\Big)^2\Big] \nonumber \\
&=\frac{4\pi\mathsf{G}}{\mathsf{c}^4}e^{-2F'+2K'}(\mathsf{c}^2\rho+3P),\label{1.52a} \\
&\frac{\partial}{\partial\varpi}\Big(\frac{e^{4F'}}{\Pi}\frac{\partial A'}{\partial\varpi}\Big)+
\frac{\partial}{\partial z}\Big(\frac{e^{4F'}}{\Pi}\frac{\partial A'}{\partial z}\Big)=0, \label{1.52b} \\
&\frac{\partial^2\Pi}{\partial\varpi^2}+\frac{\partial^2\Pi}{\partial z^2}=\frac{16\pi\mathsf{G}}{\mathsf{c}^4}
e^{-2F'+2K'}P\Pi, \label{1.52c}\\
&\frac{\partial\Pi}{\partial\varpi}\frac{\partial K'}{\partial\varpi}-
\frac{\partial\Pi}{\partial z}\frac{\partial K'}{\partial z}=
\frac{1}{2}\Big(\frac{\partial^2\Pi}{\partial \varpi^2}-
\frac{\partial^2\Pi}{\partial z^2}\Big)+
\Pi\Big[\Big(\frac{\partial F'}{\partial\varpi}\Big)^2-
\Big(\frac{\partial F'}{\partial z}\Big)^2\Big] + \nonumber \\
&-\frac{e^{4F'}}{4\Pi}
\Big[\Big(\frac{\partial A'}{\partial\varpi}\Big)^2-
\Big(\frac{\partial A'}{\partial z}\Big)^2\Big],\label{1.52d}\\
&\frac{\partial\Pi}{\partial z}\frac{\partial K'}{\partial\varpi}
+\frac{\partial\Pi}{\partial\varpi}\frac{\partial K'}{\partial z}=
\frac{\partial^2\Pi}{\partial \varpi \partial z}+2\Pi
\frac{\partial F'}{\partial\varpi}\frac{\partial F'}{\partial z}-
\frac{e^{4F'}}{2\Pi}
\frac{\partial A'}{\partial\varpi}
\frac{\partial A'}{\partial z}, \label{1.52e}\\
& \nonumber \\
&F'=-\frac{u}{\mathsf{c}^2}+\mbox{Const.}. \label{1.52f}
\end{align}
\end{subequations}

Here $\rho, P$ should be considered as given functions of $u$.

See \cite[(1.34)(1.35)(1.26)]{Meinel}.\\

Let us assume {\bf (B2)}. Then the equations \eqref{1.52d}\eqref{1.52e} can be solved as 
\begin{subequations}
\begin{align}
&\frac{\partial K'}{\partial\varpi}=
\Big[\Big(\frac{\partial\Pi}{\partial\varpi}\Big)^2+\Big(\frac{\partial\Pi}{\partial z}\Big)^2\Big]^{-1}
\Big(\frac{\partial\Pi}{\partial\varpi}\cdot\mbox{RH}\eqref{1.52d}+
\frac{\partial\Pi}{\partial z}\cdot\mbox{RH}\eqref{1.52e}\Big), \label{1.53a}\\
&\frac{\partial K'}{\partial z}=
\Big[\Big(\frac{\partial\Pi}{\partial\varpi}\Big)^2+\Big(\frac{\partial\Pi}{\partial z}\Big)^2\Big]^{-1}
\Big(-\frac{\partial\Pi}{\partial z}\cdot\mbox{RH}\eqref{1.52d}+
\frac{\partial\Pi}{\partial\varpi}\cdot\mbox{RH}\eqref{1.52e}\Big),\label{1.53b}
\end{align}
\end{subequations}
where RH\eqref{1.52d}, RH\eqref{1.52e} stand for the right-hand side of
\eqref{1.52d}, of \eqref{1.52e}, respectively.\\

Here we have a serious question: In order that $K'$ which satisfies 
\eqref{1.53a}\eqref{1.53b} exists it is necessary that the consistency condition
\begin{equation}
\frac{\partial}{\partial z}\mbox{RH}\eqref{1.53a}=\frac{\partial}{\partial\varpi}\mbox{RH}\eqref{1.53b}. \label{1.54}
\end{equation}
holds.
Is it guaranteed?
Actually
we can verify this consistency if $P=0$ by taking $\Pi=\varpi$. See \cite[\S 2.2, \S4.2]{Islam}.
But we should be careful when $P\not=0, \Pi\not=\varpi$.

By direct but tedious calculations, we can get the following observation:

\begin{Proposition}\label{Proposition1}
Let $K'$ be arbitrarily fixed and let $F', A', \Pi, \rho$ satisfy
\eqref{1.52a}, \eqref{1.52b},\eqref{1.52c} and \eqref{1.52f} with this fixed $K'$. Let us denote by
$\tilde{K}_1, \tilde{K}_3$ the right-hand sides of
\eqref{1.53a}, \eqref{1.53b}, respectively, evaluated by these $F', A', \Pi$. Then it holds that
\begin{align}
\frac{\partial\tilde{K}_1}{\partial z}-
\frac{\partial\tilde{K}_3}{\partial\varpi}&=
\frac{8\pi \mathsf{G}}{\mathsf{c}^4}
e^{-2F'+2K'}P\Pi\Big[\Big(\frac{\partial\Pi}{\partial\varpi}\Big)^2+
\Big(\frac{\partial\Pi}{\partial z}\Big)^2\Big]^{-1}\times \nonumber \\
&\times \Big[\Big(\frac{\partial K'}{\partial\varpi}-\tilde{K}_1\Big)
\frac{\partial\Pi}{\partial z}-
\Big(\frac{\partial K'}{\partial z}-\tilde{K}_3\Big)
\frac{\partial\Pi}{\partial\varpi}\Big].\label{1.55}
\end{align}
\end{Proposition}

Therefore, as a conclusion, if $K'$ satisfies \eqref{1.53a} \eqref{1.53b},
then the consistency condition \eqref{1.54} holds, since 
$$\frac{\partial K'}{\partial\varpi}=\tilde{K}_1,\qquad
\frac{\partial K'}{\partial z}=\tilde{K}_3.
$$ Of course this conclusion in itself is a vicious circular argument of no use. However the following argument is useful:

\begin{Lemma}\label{Lemma1}
Let us consider a bounded disk $\mathfrak{D}=\{
(\varpi, z)\  |\  r=\sqrt{\varpi^2+z^2}<R \}$ of $(\varpi,z)$-plane and denote by $\bar{\mathfrak{D}}$ the closure of $\mathfrak{D}$. Suppose that $K' \in C(\bar{\mathfrak{D}})$ is given and that $F', A',  \rho \in C^2(\bar{\mathfrak{D}}), 
\Pi \in C^3(\bar{\mathfrak{D}})$ satisfy \eqref{1.52a},\eqref{1.52b},\eqref{1.52c} and \eqref{1.52f} with this $K'$. Suppose {\bf (B2)} holds on
$\bar{\mathfrak{D}}$. Let us denote by
$\tilde{K}_1,\tilde{K}_3$ the right-hand sides of
\eqref{1.53a},\eqref{1.53b}, respectively, evaluated by these
$F', A', \Pi$. (They are $C^1$-functions on $\bar{\mathfrak{D}}$.)
Put
\begin{equation}
\tilde{K}(\varpi, z):=
\int_0^z\tilde{K}_3(0,z')dz'+
\int_0^{\varpi}\tilde{K}_1(\varpi',z)d\varpi'\label{1.56}
\end{equation}
for $(\varpi, z) \in \mathfrak{D}$. If $\tilde{K}=K'$, then $K'$
satisfies 
\begin{equation}
\frac{\partial K'}{\partial\varpi}=\tilde{K}_1,\qquad
\frac{\partial K'}{\partial z}=\tilde{K}_3,\label{1.57}
\end{equation}
that is, the equations \eqref{1.52d}\eqref{1.52e} are satisfied.
\end{Lemma}

{\bf Proof.}  Suppose that $\tilde{K}=K'$. It follows from \eqref{1.56} with $\tilde{K}=K'$ that
\begin{align}
\frac{\partial K'}{\partial \varpi}(\varpi, z)&=\tilde{K}_1(\varpi, z), \label{1.58} \\
\frac{\partial K'}{\partial z}(\varpi, z)&=
\tilde{K}_3(0,z)+\int_0^{\varpi}
\frac{\partial\tilde{K}_1}{\partial z}(\varpi', z)d\varpi'.\label{1.59}
\end{align}
Put
\begin{equation}
L(\varpi, z):=\frac{\partial \tilde{K}_1}{\partial z}-
\frac{\partial\tilde{K}_3}{\partial\varpi},
\end{equation}
which is a continuous function on $\bar{\mathfrak{D}}$. Then \eqref{1.59} reads
\begin{equation}
\frac{\partial K'}{\partial z}(\varpi, z)=
\tilde{K}_3(\varpi, z)+
\int_0^{\varpi}
L(\varpi', z)d\varpi'.\label{1.61}
\end{equation}
Now therefore \eqref{1.55} of Proposition~\ref{Proposition1} reads
\begin{equation}
L(\varpi, z)=-
\frac{8\pi\mathsf{G}}{\mathsf{c}^4}
e^{-2F'+2K'}P\Pi\Big[\Big(\frac{\partial\Pi}{\partial\varpi}\Big)^2+
\Big(\frac{\partial\Pi}{\partial z}\Big)^2\Big]^{-1}\frac{\partial\Pi}{\partial\varpi}
\int_0^{\varpi}L(\varpi',z)d\varpi'.
\end{equation}
Since the function
$$
\frac{8\pi\mathsf{G}}{\mathsf{c}^4}
e^{-2F'+2K'}P\Pi\Big[\Big(\frac{\partial\Pi}{\partial\varpi}\Big)^2+
\Big(\frac{\partial\Pi}{\partial z}\Big)^2\Big]^{-1}\frac{\partial\Pi}{\partial\varpi}
$$
is bounded on the compact $\bar{\mathfrak{D}}$ thanks to {\bf (B2)}, the Gronwall's argument implies that
$L(\varpi, z)=0$ on $\mathfrak{D}$ so that \eqref{1.61} reads
\begin{equation}
\frac{\partial K'}{\partial z}(\varpi, z)=\tilde{K}_3(\varpi, z).\label{1.63}
\end{equation}
Thus \eqref{1.58} and \eqref{1.63} complete the proof. $\square$\\



\section{Post-Newtonian approximation}

We are going to find a solution of 
\eqref{1.52a}\eqref{1.52b}\eqref{1.52c}\eqref{1.52d}\eqref{1.52e}\eqref{1.52f} which approaches to the solution of the Euler-Poisson equations constructed in \cite{JJTM} as $\mathsf{c} \rightarrow +\infty$.\\


Recall that we are assuming 
{\bf (A)}.


Then there are smooth functions $\Lambda_{u}, \Lambda_{\rho},\Lambda_P$ which vanish at 0 such that
\begin{subequations}
\begin{align}
&u:=\int_0^{\rho}\frac{dP}{\rho+P/\mathsf{c}^2}=
\frac{\mathsf{A}\gamma}{\gamma-1}\rho^{\gamma-1}
(1+\Lambda_u(\mathsf{A}\rho^{\gamma-1}/\mathsf{c}^2))\quad\mbox{for}\quad \rho >0, \\
&\rho=\Big(\frac{\gamma-1}{\mathsf{A}\gamma}\Big)^{\frac{1}{\gamma-1}}
u^{\frac{1}{\gamma-1}}(1+\Lambda_{\rho}(u/\mathsf{c}^2))\quad\mbox{for}\quad
u>0, \\
&P=\mathsf{A}^{-\frac{1}{\gamma-1}}\Big(\frac{\gamma-1}{\gamma}\Big)^{\frac{\gamma}{\gamma-1}}
u^{\frac{\gamma}{\gamma-1}}(1+
\Lambda_P(u/\mathsf{c}^2))\quad \mbox{for}\quad u>0.
\end{align}
\end{subequations}

Now let us fix the Newtonian limit

\begin{subequations}
\begin{align}
&u_{\mathsf{N}}=u_{\mathsf{O}}\Theta\Big(\frac{r}{\mathsf{a}},\zeta;\frac{1}{\gamma-1}, \mathsf{b}\Big),\\
&\mathsf{a}=\frac{1}{\sqrt{4\pi\mathsf{G}}}
\Big(\frac{\mathsf{A}\gamma}{\gamma-1}\Big)^{\frac{1}{2(\gamma-1)}}
u_{\mathsf{O}}^{-\frac{2-\gamma}{2(\gamma-1)}},
\quad
\mathsf{b}=
\frac{\Omega^2}{2\pi\mathsf{G}}
\Big(\frac{\mathsf{A}\gamma}{\gamma-1}\Big)^{\frac{1}{\gamma-1}}
u_{\mathsf{O}}^{-\frac{1}{\gamma-1}},
\end{align}
\end{subequations}
where $u_{\mathsf{O}}=u_{\mathsf{N}}(\vec{0})$ is a given positive number, 
$$r=\sqrt{\varpi^2+z^2},\qquad \zeta=\frac{z}{\sqrt{\varpi^2+z^2}},
$$
and $\Theta$ is the distorted Lane-Emden function of index $\frac{1}{\gamma-1}$ with parameter $\mathsf{b}$ which is supposed to be sufficiently small. See \cite{JJTM} and \cite{JJTM.Ext}. Therefore we consider it on the domain
$$\mathfrak{D}=\{ (\varpi, z) \  |\  r=\sqrt{\varpi^2+z^2}< R \}, $$
in which the support of $u_{\mathsf{N}}$ is included. Here $R$ is fixed so that $ 2\mathsf{a}\xi_1(\frac{1}{\gamma-1})\leq R\leq \mathsf{a}\Xi_0$. 

So, we consider the perturbation $w$ defined by
\begin{equation}
u=u_{\mathsf{N}}+\frac{w}{\mathsf{c}^2},
\end{equation}
$w$ being an unknown function defined on $\mathfrak{D}$
which satisfies
\begin{equation}
w(0,0)=0.
\end{equation} Then
\begin{subequations}
\begin{align}
&\rho=\Big(\frac{\gamma-1}{\mathsf{A}\gamma}\Big)^{\frac{1}{\gamma-1}}
(u\vee 0)^{\frac{1}{\gamma-1}}(1+\Lambda_{\rho}(u/\mathsf{c}^2)), \\
&P=\mathsf{A}^{-\frac{1}{\gamma-1}}
\Big(\frac{\gamma-1}{\gamma}\Big)^{\frac{\gamma}{\gamma-1}}
(u\vee 0)^{\frac{\gamma}{\gamma-1}}(1+
\Lambda_P(u/\mathsf{c}^2)),
\end{align}
\end{subequations}
where
\begin{equation}
u \vee 0=\max\{ u, 0\}=\max\{ u_{\mathsf{N}}+\frac{w}{\mathsf{c}^2}, 0\}.
\end{equation}
Of course the Newtonian limits $\rho_{\mathsf{N}}, P_{\mathsf{N}}$ are
\begin{subequations}
\begin{align}
&\rho_{\mathsf{N}}=\Big(\frac{\gamma-1}{A\gamma}\Big)^{\frac{1}{\gamma-1}}
(u_{\mathsf{N}}\vee 0)^{\frac{1}{\gamma-1}}, \\
&P_{\mathsf{N}}=A^{-\frac{1}{\gamma-1}}
\Big(\frac{\gamma-1}{\gamma}\Big)^{\frac{\gamma}{\gamma-1}}
(u_{\mathsf{N}}\vee 0)^{\frac{\gamma}{\gamma-1}}.
\end{align}
\end{subequations}\\

Moreover we see
\begin{align}
\mathsf{c}^2(\rho-\rho_{\mathsf{N}})&=\frac{1}{\gamma-1}\frac{\rho_{\mathsf{N}}}{u_{\mathsf{N}}}w + \lambda_1\rho_{\mathsf{N}}u_{\mathsf{N}} + \nonumber \\
&+\mathsf{c}^2H_{\rho}(w)+\frac{\lambda_1}{\mathsf{c}^2}\Big[\rho_{\mathsf{N}}w+
\frac{1}{\gamma-1}\frac{\rho_{\mathsf{N}}}{u_{\mathsf{N}}}wu+\mathsf{c}^2H_{\rho}(w)u\Big]+ \nonumber \\
&+\frac{\lambda_2}{\mathsf{c}^2}f_{\rho}(u)u^2\Big(1+\Lambda''(u/\mathsf{c}^2)\Big),
\end{align}
where 
$$u=u_{\mathsf{N}}+\frac{w}{\mathsf{c}^2} $$
and the function $H_{\rho}(w)$ is defined by
\begin{equation}
H_{\rho}(w)=f_{\rho}\Big(u_{\mathsf{N}}+\frac{w}{\mathsf{c}^2}\Big)-f_{\rho}(u_{\mathsf{N}})-Df_{\rho}(u_{\mathsf{N}})\frac{w}{\mathsf{c}^2}
\end{equation}
from the function $f_{\rho}(u)$, which gives $\rho_{\mathsf{N}}$ from
$u_{\mathsf{N}}$, :
\begin{equation}
f_{\rho}(u)=
\Big(\frac{\gamma-1}{\mathsf{A}\gamma}\Big)^{\frac{1}{\gamma-1}}
(u\vee 0)^{\frac{1}{\gamma-1}}.
\end{equation}
The constants $\lambda_1, \lambda_2$ are those appearing in the expression of $\Lambda_{\rho}$:
\begin{equation}
\Lambda_{\rho}(\xi)=\lambda_1\xi+\lambda_2\xi^2(1+\Lambda''(\xi)), \label{3.11}
\end{equation}
where $\Lambda''(0)=0$.
Here we note that
$$Df_{\rho}(u_{\mathsf{N}})=\frac{1}{\gamma-1}\frac{\rho_{\mathsf{N}}}{u_{\mathsf{N}}}.
$$\\

Let us introduce the variables $V, X, Y$ to put
\begin{align}
&\mathsf{c}^2F'=\Phi_{\mathsf{N}}-\frac{\Omega^2}{2}\varpi^2-\frac{w}{\mathsf{c}^2}, \\
&\mathsf{c}^2K'=-\frac{\Omega^2}{2}\varpi^2+\frac{V}{\mathsf{c}^2}, \\
&\Pi=\varpi\Big(1+\frac{X}{\mathsf{c}^4}\Big), \\
&A'=-\frac{\Omega}{\mathsf{c}}\varpi^2\Big(1+\frac{Y}{\mathsf{c}^2}\Big).
\end{align}

Here let us recall
\begin{equation}
u_{\mathsf{N}}+\Phi_{\mathsf{N}}-\frac{\Omega^2}{2}\varpi^2=\mbox{Const.}
\end{equation}
We shall denote
\begin{equation}
\Phi_{\mathsf{N}}':=\Phi_{\mathsf{N}}-\frac{\Omega^2}{2}\varpi^2,
\end{equation}
and
\begin{equation}
\Phi':=\mathsf{c}^2F'=\Phi_{\mathsf{N}}'-\frac{w}{\mathsf{c}^2}
=\Phi_{\mathsf{N}}-\frac{\Omega^2}{2}\varpi^2-\frac{w}{\mathsf{c}^2}.
\end{equation}



We are going to rewrite the equations
for the unknown variables $w, V, X, Y$.\\

The equations \eqref{1.52a}\eqref{1.52b}\eqref{1.52c} are reduced to
\begin{subequations}
\begin{align}
&\frac{\partial^2w}{\partial\varpi^2}+\frac{1}{\varpi}\frac{\partial w}{\partial\varpi}+
\frac{\partial^2w}{\partial z^2}
+\frac{4\pi\mathsf{G}}{\gamma-1}\frac{\rho_{\mathsf{N}}}{u_{\mathsf{N}}}w
-8\Phi_{\mathsf{N}}'\Omega^2-2\Omega^2 Y_1+ \nonumber \\
&+4\pi\mathsf{G}\Big(-2\Phi_{\mathsf{N}}\rho_{\mathsf{N}}+
\lambda_1\rho_{\mathsf{N}}u_{\mathsf{N}}+3P) + R_a=0
\quad\mbox{with}\quad w(0,0)=0, \label{Eq.w}\\
&\frac{\partial^2Y}{\partial\varpi^2}+\frac{3}{\varpi}\frac{\partial Y}{\partial\varpi}+
\frac{\partial^2Y}{\partial z^2}+
\frac{8}{\varpi}\frac{\partial\Phi_{\mathsf{N}}'}{\partial\varpi}+R_b=0, \label{Eq.Y}\\
&\frac{\partial^2X}{\partial\varpi^2}+\frac{2}{\varpi}\frac{\partial X}{\partial\varpi}+
\frac{\partial^2X}{\partial z^2}-16\pi\mathsf{G}P_{\mathsf{N}}+R_c=0,
\label{Eq.X}
\end{align}
\end{subequations}
where 
$\lambda_1$ is the constant in \eqref{3.11}.
The terms $R_a, R_b, R_c$ are as following:

\begin{align}
R_a&:=-\Big(1+\frac{X}{\mathsf{c}^4}\Big)^{-1}\frac{1}{\mathsf{c}^2}\varpi\frac{\partial X}{\partial\varpi}
\frac{1}{\varpi}\frac{\partial\Phi'}{\partial\varpi} \nonumber \\
&-\Big[\Big(1+\frac{X}{\mathsf{c}^4}\Big)^{-1}e^{4F'}
-1-\frac{4\Phi_{\mathsf{N}}'}{\mathsf{c}^2}\Big]\cdot 2\mathsf{c}^2\Omega^2 \nonumber \\
&-\Big[\Big(1+\frac{X}{\mathsf{c}^4}\Big)^{-1}e^{4F'}
-1\Big]\cdot 2\Omega^2{Y_1} \nonumber \\
&-\Big(1+\frac{X}{\mathsf{c}^4}\Big)^{-1}e^{4F'}
\cdot \frac{\Omega^2}{2\mathsf{c}^2}\Big[(Y_1)^2+\varpi^2(Y_3)^2\Big], \nonumber \\
&+4\pi\mathsf{G}\Big[\mathsf{c}^2(\rho-\rho_{\mathsf{N}})-\frac{1}{\gamma-1}
\frac{\rho_{\mathsf{N}}}{u_{\mathsf{N}}}w-\lambda_1\rho_{\mathsf{N}}u_{\mathsf{N}} \nonumber \\
&+\mathsf{c}^2(e^{2(K'-F')}-1)(\rho-\rho_{\mathsf{N}})
+\mathsf{c}^2\Big(e^{2(K'-F')}-1+
\frac{2\Phi_{\mathsf{N}}}{\mathsf{c}^2}\Big)\rho_{\mathsf{N}} \nonumber \\
&+3(e^{2(K'-F')}-1)P+3(P-P_{\mathsf{N}})\Big].
\label{2.21a}
\end{align}
\begin{align}
R_b&:=-\frac{8}{\mathsf{c}^2}\frac{1}{\varpi}\frac{\partial w}{\partial \varpi} 
-\Big(1+\frac{X}{\mathsf{c}^4}\Big)^{-1}
\frac{1}{\mathsf{c}^2}\frac{1}{\varpi}\frac{\partial X}{\partial\varpi}\Big(2+\frac{Y_1}{\mathsf{c}^2}\Big)
+\frac{1}{\mathsf{c}^4}\Big(1+\frac{X}{\mathsf{c}^4}\Big)^{-1}
\frac{\partial X}{\partial z}Y_3 \nonumber \\
&+\frac{4}{\mathsf{c}^2}\frac{1}{\varpi}\frac{\partial\Phi'}{\partial\varpi}
\cdot Y_1
+\frac{4}{\mathsf{c}^2}\frac{\partial\Phi'}{\partial z}\cdot Y_3, \label{2.21b}\\
R_c&:=16\pi\mathsf{G}\Big[
-e^{2(K'-F')}P+P_{\mathsf{N}}
-\frac{1}{\mathsf{c}^4}e^{2(K'-F')}PX\Big].\label{2.21c}
\end{align}

Here we should read 
\begin{align}
&\mathsf{c}^2F'=\Phi'=\Phi_N'-\frac{w}{\mathsf{c}^2}=
-\frac{\Omega^2}{2}\varpi^2+\Phi_{\mathsf{N}}-\frac{w}{\mathsf{c}^2}, \label{PH}\\
&\mathsf{c}^2(K'-F')=-\frac{\Omega^2}{2}\varpi^2
-\Phi_{\mathsf{N}}'+\frac{1}{\mathsf{c}^2}(V+w)
=-\Phi_{\mathsf{N}}+\frac{1}{\mathsf{c}^2}(V+w). \label{PS}
\end{align}
Note that  we have used the identity
\begin{equation}
\Big(1+\frac{X}{\mathsf{c}^4}\Big)^{-1}\Big(1+\frac{X_1}{\mathsf{c}^4}\Big)-1=
\Big(1+\frac{X}{\mathsf{c}^4}\Big)^{-1}\frac{1}{\mathsf{c}^4}
\varpi\frac{\partial X}{\partial\varpi},
\end{equation}
while $$X_1:=X+\varpi\frac{\partial X}{\partial\varpi} $$
gives
$$\frac{\partial\Pi}{\partial\varpi}=1+\frac{X_1}{\mathsf{c}^4}.$$
We have denoted 
\begin{equation}
Y_1=2Y+\varpi\frac{\partial Y}{\partial\varpi},\qquad Y_3=\frac{\partial Y}{\partial z}.
\end{equation}\\

Now
the equations for $V$ are
\begin{subequations}
\begin{align}
&\frac{\partial V}{\partial\varpi}=\mathsf{c}^2\Omega^2\varpi+
\mathsf{c}^4\cdot\mbox{RH}\eqref{1.53a}, \label{2.25a}\\
&\frac{\partial V}{\partial z}=\mathsf{c}^4\cdot\mbox{RH}\eqref{1.53b}.
\label{2.25b}
\end{align}
\end{subequations}
These read
\begin{subequations}
\begin{align}
\frac{\partial V}{\partial\varpi}&=-4\Phi_{\mathsf{N}}'\varpi\Omega^2+
\frac{\partial X}{\partial\varpi}+\frac{\varpi}{2}\Big(
\frac{\partial^2 X}{\partial\varpi^2}-\frac{\partial^2 X}{\partial z^2}\Big)+
\nonumber \\
 &+\varpi\Big(\Big(\frac{\partial\Phi_{\mathsf{N}}'}{\partial\varpi}\Big)^2
-\Big(\frac{\partial\Phi_{\mathsf{N}}'}{\partial z}\Big)^2\Big) +\varpi\Omega Y_1+R_d, \label{Eq.K1}\\
\frac{\partial V}{\partial z}&=
\frac{\partial X}{\partial z}+\varpi\frac{\partial^2X}{\partial\varpi\partial z}+
2\varpi\frac{\partial\Phi_{\mathsf{N}}'}{\partial\varpi}\frac{\partial\Phi_{\mathsf{N}}'}{\partial z}+
\varpi^2\Omega Y_3+R_e, \label{Eq.K3}
\end{align}
\end{subequations}
where
\begin{align}
R_d&=
\mathsf{c}^2\Omega^2\varpi\Big[1-
\Big(1+\frac{X_*}{\mathsf{c}^4}\Big)^{-1}\Big(1+\frac{X_1}{\mathsf{c}^4}\Big)
\Big(1+\frac{X}{\mathsf{c}^4}\Big)^{-1}e^{4F'}+\frac{4\Phi_N'}{\mathsf{c}^2}\Big] \nonumber \\
&+\frac{1}{2}\Big[\Big(1+\frac{X_*}{\mathsf{c}^4}\Big)^{-1}\Big(1+
\frac{X_1}{\mathsf{c}^4}\Big)-1\Big]
\Big(2\frac{\partial X}{\partial\varpi}+\varpi\frac{\partial^2X}{\partial\varpi^2}-
\varpi\frac{\partial^2X}{\partial z^2}\Big) \nonumber \\
&+\Big[\Big(1+\frac{X_*}{\mathsf{c}^4}\Big)^{-1}\Big(1+\frac{X_1}{\mathsf{c}^4}\Big)
\Big(1+\frac{X}{\mathsf{c}^4}\Big)-1\Big]\varpi\Big[
\Big(\frac{\partial\Phi_N'}{\partial\varpi}\Big)^2
-\Big(\frac{\partial\Phi_N'}{\partial z}\Big)^2\Big] \nonumber \\
&-\Big[\Big(1+\frac{X_*}{\mathsf{c}^4}\Big)^{-1}\Big(1+\frac{X_1}{\mathsf{c}^4}\Big)
\Big(1+\frac{X_1}{\mathsf{c}^4}\Big)e^{4F'}-1\Big]\varpi\Omega^2 Y_1 \nonumber \\
&+\Big(1+\frac{X_*}{\mathsf{c}^4}\Big)^{-1}\frac{\varpi}{\mathsf{c}^2}\frac{\partial X}{\partial z}
\Big[\frac{1}{\mathsf{c}^2}\Big(\frac{\partial X}{\partial z}+\varpi
\frac{\partial^2X}{\partial\varpi\partial z}\Big) \nonumber \\
&+\frac{2\varpi}{\mathsf{c}^2}\Big(1+\frac{X}{\mathsf{c}^4}\Big)
\frac{\partial\Phi'}{\partial\varpi}\frac{\partial\Phi'}{\partial z} \nonumber \\
&
-\frac{e^{4F'}}{2}\Big(1+\frac{X}{\mathsf{c}^4}\Big)^{-1}\frac{\varpi^2}{\mathsf{c}^2}\Big(-2\Omega+
\frac{\Omega^2}{\mathsf{c}^2}Y_1Y_3\Big)\Big],
\label{Rd}
\\
R_e&=
\Big[\Big(1+\frac{X_*}{\mathsf{c}^4}\Big)^{-1}\Big(1+\frac{X_1}{\mathsf{c}^4}\Big)-1\Big]
\Big(\frac{\partial X}{\partial\varpi}+\varpi\frac{\partial^2X}{\partial\varpi\partial z}\Big) \nonumber \\
&+\Big[\Big(1+\frac{X_*}{\mathsf{c}^4}\Big)^{-1}\Big(1+\frac{X_1}{\mathsf{c}^4}\Big)
\Big(1+\frac{X}{\mathsf{c}^4}\Big)-1\Big]2\varpi
\frac{\partial\Phi_N'}{\partial\varpi}\frac{\partial\Phi_N'}{\partial z} \nonumber \\
&-\Big[\Big(1+\frac{X_*}{\mathsf{c}^4}\Big)^{-1}\Big(1+\frac{X_1}{\mathsf{c}^4}\Big)
\Big(1+\frac{X}{\mathsf{c}^4}\Big)^{-1}e^{4F'}-1\Big]
\Omega^2\varpi^2 Y_3 \nonumber \\
&-\Big(1+\frac{X_*}{\mathsf{c}^4}\Big)^{-1}
\frac{\varpi}{\mathsf{c}^2}
\frac{\partial X}{\partial z}\Big[
\frac{1}{2\mathsf{c}^2}\Big(2\frac{\partial X}{\partial\varpi}
+\varpi\frac{\partial^2X}{\partial\varpi^2}-
\varpi\frac{\partial^2X}{\partial z^2}\Big) \nonumber\\
&+\frac{\varpi}{\mathsf{c}^2}\Big(1+\frac{X}{\mathsf{c}^4}\Big)\frac{\partial\Phi'}{\partial\varpi}
\frac{\partial\Phi'}{\partial z} \nonumber \\
&-\frac{e^{4F'}}{4}
\Big(1+\frac{X}{\mathsf{c}^4}\Big)^{-1}\varpi\Big(
\Big(\Omega^2\Big(2+\frac{Y_1}{\mathsf{c}^2}\Big)^2-
\frac{\varpi^2}{\mathsf{c}^4}\Omega^2Y_3^2\Big)\Big]. \label{Re}
\end{align}
Here we put
\begin{equation}
X_*:=2\Big(X+\varpi\frac{\partial X}{\partial\varpi}\Big)+
\frac{1}{\mathsf{c}^4}
\Big[X^2+2X\varpi\frac{\partial X}{\partial\varpi}+
\varpi^2\Big(\Big(\frac{\partial X}{\partial\varpi}\Big)^2+
\Big(\frac{\partial X}{\partial z}\Big)^2\Big)\Big]
\end{equation}
so that
\begin{equation}
\Big(\frac{\partial\Pi}{\partial\varpi}\Big)^2+\Big(\frac{\partial\Pi}{\partial z}\Big)^2=1+\frac{X_*}{\mathsf{c}^4}.
\end{equation}\\

Of course, we can restate Lemma~\ref{Lemma1}
as
\begin{Lemma}\label{Lemma2}
Let $\mathfrak{D}, \bar{\mathfrak{D}}$ be those of Lemma~\ref{Lemma1}.
Suppose that $V \in C(\bar{\mathfrak{D}})$ is given and that
$ Y,  w \in C^2(\bar{\mathfrak{D}}),
X \in C^3(\bar{\mathfrak{D}})$ satisfy
\eqref{Eq.w},\eqref{Eq.Y},\eqref{Eq.X}  with this $V$.
Here, of course, the terms $\rho, P, K'-F'$ in \eqref{Eq.w},\eqref{Eq.X} should read
\begin{align*}
&\rho=\Big(\frac{\gamma-1}{A\gamma}\Big)^{\frac{1}{\gamma-1}}
(u\vee 0)^{\frac{1}{\gamma-1}}(1+\Lambda_{\rho}(u/\mathsf{c}^2)), \\
&P=A^{-\frac{1}{\gamma-1}}
\Big(\frac{\gamma-1}{\gamma}\Big)^{\frac{\gamma}{\gamma-1}}
(u\vee 0)^{\frac{\gamma}{\gamma-1}}(1+
\Lambda_P(u/\mathsf{c}^2)),\\
\mbox{with}&\quad u=u_{\mathsf{N}}+\frac{w}{\mathsf{c}^2}, \\
&\mathsf{c}^2(K'-F')=-\Phi_{\mathsf{N}}+\frac{1}{\mathsf{c}^2}(V+w).
\end{align*}
Let us denote by $\tilde{V}_1, \tilde{V}_3$ the right-hand sides of
\eqref{Eq.K1},\eqref{Eq.K3}, respectively, evaluated by these
$w, Y, X$. Put
\begin{equation}
\tilde{V}:=\int_0^z\tilde{V}_3(0,z')dz'+
\int_0^{\varpi}\tilde{V}_1(\varpi',z)d\varpi'.
\end{equation}
If $\tilde{V}=V$, then $V$ satisfies
\begin{equation}
\frac{\partial V}{\partial\varpi}=\tilde{V}_1,
\qquad
\frac{\partial V}{\partial z}=\tilde{V}_3,
\end{equation}
that is, the equations \eqref{Eq.K1},\eqref{Eq.K3} are satisfied.
\end{Lemma}

\section{Functional spaces}

We are going to prepare some notations to prove the existence of solutions.

\subsection{H\"{o}lder spaces and Newtonian potentials}

Let $n=3,4,5$. For given $\Xi>0$, we denote
$$ B^{(n)}(\Xi)=\{ \xi=(\xi_1,\cdots,\xi_n) \in \mathbb{R}^n\  |\  |\xi|=\sqrt{\sum_j(\xi_j)^2}<\Xi \} $$
and $\bar{B}^{(n)}(\Xi)=\{|\xi|\leq\Xi\} $.

For a continuous function $f$ on $\bar{B}^{(n)}(\Xi)$ and $l=0,1,2$, we put
$$\| f ; C^l(\bar{B}^{(n)}(\Xi))\|:=\sum_{|L|\leq l}\sup_{|\xi|\leq\Xi}|\partial_{\xi}^Lf(\xi)|,$$
where
$$ \partial_{\xi}^L =\Big(\frac{\partial}{\partial\xi_1}\Big)^{L_1}\cdots
\Big(\frac{\partial}{\partial \xi_n}\Big)^{L_n} $$
for $L=(L_1,\cdots, L_n)$ and $|L|=L_1+\cdots +L_n$.
Then $\| \cdot ; C^l(\bar{B}^{(n)}(\Xi))\|$ is the usual norm of the Banach space
$C^l(\bar{B}^{(n)}(\Xi))$.\\

Let us fix a number $\alpha$ such that
\begin{equation}
0<\alpha <\min\{\frac{1}{\gamma-1}-1, 1\}. 
\end{equation}

We put
\begin{align*}
\|f ; C^{l,\alpha}(\bar{B}^{(n)}(\Xi))\|:=&
\|f;C^l(\bar{B}^{(n)}(\Xi))\|+ \nonumber \\
&+\sup_{|\xi'|,|\xi|\leq \Xi, |L|=l}
\frac{|\partial_{\xi}^Lf(\xi')-\partial_{\xi}^Lf(\xi)|}{|\xi'-\xi|^{\alpha}}
\end{align*}
and
$$C^{l,\alpha}(\bar{B}^{(n)}(\Xi))=\{ f \in C^l(\bar{B}^{(n)}(\Xi))\  |\  
\|f; C^{l,\alpha}(\bar{B}^{(n)}(\Xi))\|<\infty \}.
$$

Then we have the following 
\begin{Proposition}\label{Prop6}
i) The bilinear mapping $(f,g)\mapsto f\cdot g$ is continuous as
$$C^l(\bar{B}^{(n)}(\Xi))\times C^{l'}(\bar{B}^{(n)}(\Xi)) \rightarrow
C^{l\wedge l'}(\bar{B}^{(n)}(\Xi))$$
and as
$$C^{l,\alpha}(\bar{B}^{(n)}(\Xi))\times C^{l,\alpha}(\bar{B}^{(n)}(\Xi)) \rightarrow
C^{l,\alpha}(\bar{B}^{(n)}(\Xi))$$

ii) If $l\leq l'$, then the imbedding 
$C^{l'}(\bar{B}^{(n)}(\Xi))\rightarrow C^l(\bar{B}^{(n)}(\Xi))$ and
$C^{l',\alpha}(\bar{B}^{(n)}(\Xi))\rightarrow C^{l,\alpha}(\bar{B}^{(n)}(\Xi))$
are continuous.

iii) The imbedding $C^{l+1}(\bar{B}^{(n)}(\Xi))\rightarrow
C^{l,\alpha}(\bar{B}^{(n)}(\Xi))$ is continuous.

iv) The imbedding $C^{l,\alpha}(\bar{B}^{(n)}(\Xi)) \rightarrow C^l
(\bar{B}^{(n)}(\Xi))$ is compact.
\end{Proposition}

Let us fix $\Xi$ and $\Xi_0$ such that 
$$2\xi_1\Big(\frac{1}{\gamma-1}\Big)\leq \Xi <\Xi_0$$
and a cut off function $\chi\in C^{\infty}([0,+\infty[)$ such that
$\chi(\eta)=1$ for $0\leq\eta\leq \Xi$, $0<\chi(\eta)<1$ for $\Xi<\eta<(\Xi+\Xi_0)/2$, and
$\chi(\eta)=0$ for $(\Xi+\Xi_0)/2\leq \eta$. 

For $g \in C^0(\bar{B}^{(n)}(\Xi_0))$ we put
\begin{equation}
\mathcal{K}^{(n)}g(\xi)=
\frac{1}{S_n}\int\frac{g(\xi')\chi(|\xi'|)}{|\xi-\xi'|^{n-2}}d\xi',
\end{equation}
with $S_n=2(n-2)\pi^{n/2}/\Gamma(n/2))$. Then we have the following

\begin{Proposition}\label{Prop.FS3}
i) The linear operator $\mathcal{K}^{(n)}$ is continuous as
$$ C^0(\bar{B}^{(n)}(\Xi_0))\rightarrow C^1(\bar{B}^{(n)}(\Xi_0))$$
and
$$ C^{0,\alpha}(\bar{B}^{(n)}(\Xi_0))\rightarrow
C^{2,\alpha}(\bar{B}^{(n)}(\Xi_0))$$.

ii) If $g \in C^{0,\alpha}(\bar{B}^{(n)}(\Xi_0))$, then
$f=\mathcal{K}^{(n)}g \in C^{2,\alpha}(\bar{B}^{(n)}(\Xi_0))$ satisfies 
$$\triangle f + g=0$$
on $B^{(n)}(\Xi))$. Here $\triangle$ denotes the $n$-dimensional  Laplace operator
$$\triangle=\triangle^{(n)}=\sum_{j=1}^n
\Big(\frac{\partial}{\partial\xi_j}\Big)^2.
$$
\end{Proposition}
See , e.g., \cite[Theorem 4.5]{GilbergT}.

\subsection{Dimensionless functional spaces and solution of the Poisson equations}

Let us use the positive parameter $\mathsf{a}$. 

For a function $Q$ of $(\varpi, x)$ and $n=3,4,5$, we consider the function $Q^{\flat(n)}$ of $\xi \in \mathbb{R}^n$ defind by
$$Q^{\flat(n)}(\xi)=Q(\varpi, z)$$
with
$$\xi=\frac{x}{\mathsf{a}},\qquad \varpi=\sqrt{(x_1)^2+\cdots(x_{n-1})^2},\quad
z=x_n.
$$

Considering functions on the region
$$\bar{\mathfrak{D}}(R)):=\{(\varpi, x)\  |\  r=\sqrt{\varpi^2+z^2}\leq R \}, $$
we put
\begin{align*}
&\mathfrak{C}^l(R):=\{ Q \in C(\bar{\mathfrak{D}}(R)) \  |\  
Q(\varpi, -z)=Q(\varpi, z)\quad\forall z \\
&\mbox{and}\quad Q^{\flat(n)}\in C^l(\bar{B}^{(n)}(R/\mathsf{a})) \}, \\
&\|Q; \mathfrak{C}^l(R)\|:=\|Q^{\flat(n)}; C^l(\bar{B}^{(n)}(R/\mathsf{a}))\|, \\
&\mathfrak{C}^{l,\alpha}(R):=\{ Q \in C(\bar{\mathfrak{D}}(R)) \  |\  
Q(\varpi, -z)=Q(\varpi, z)\quad\forall z \\
&\mbox{and}\quad Q^{\flat(n)}\in C^{l,\alpha}(\bar{B}^{(n)}(R/\mathsf{a})) \}, \\
&\|Q; \mathfrak{C}^{l,\alpha}(R)\|:=\|Q^{\flat(n)}; C^{l,\alpha}(\bar{B}^{(n)}(R/\mathsf{a}))\|
\end{align*}

It is easy to see that the spaces $\mathfrak{C}^l(R),
\mathfrak{C}^{l,\alpha}(R)$, and the norms
$\|\cdot ; \mathfrak{C}^l(R)\|, 
\|\cdot; \mathfrak{C}^{l,\alpha}\|$ do not depend on the choice of $n$, but depend only on $\mathsf{a}$.

Proposition~\ref{Prop6} reads
\begin{Proposition}
i) The bilinear mapping $(f,g)\mapsto f\cdot g$ is continuous as
$$\mathfrak{C}^l(R)\times \mathfrak{C}^{l'}(R) \rightarrow
\mathfrak{C}^{l\wedge l'}(R))$$
and as
$$\mathfrak{C}^{l,\alpha}(R)\times \mathfrak{C}^{l,\alpha}(R) \rightarrow
\mathfrak{C}^{l,\alpha}(R)$$

ii) If $l\leq l'$, then the imbedding 
$\mathfrak{C}^{l'}(R)\rightarrow \mathfrak{C}^l(R)$ and
$\mathfrak{C}^{l',\alpha}(R)\rightarrow \mathfrak{C}^{l,\alpha}(R)$
are continuous.

iii) The imbedding $\mathfrak{C}^{l+1}(R)\rightarrow
\mathfrak{C}^{l,\alpha}(R)$ is continuous.

iv) The imbedding $\mathfrak{C}^{l,\alpha}(R) \rightarrow \mathfrak{C}^l
(R)$ is compact.
\end{Proposition}


Let us fix $R$ such that
$$2\xi_1\Big(\frac{1}{\gamma-1}\Big)\leq \Xi=\frac{R}{\mathsf{a}} <\Xi_0 $$
and put $R_0=\mathsf{a} \Xi_0$.

We define the operator $\mathfrak{K}^{(n)}$ by
$$(\mathfrak{K}^{(n)}Q)^{\flat(n)}=\mathcal{K}^{(n)}Q^{\flat(n)}. $$
The operator is well-defined, and turns out to be a compact linear operator in $\mathfrak{C}^0(R_0)$. Proposition~\ref{Prop.FS3} reads

\begin{Proposition}
i) The operator $\mathfrak{K}^{(n)}$ is continuous from
$\mathfrak{C}^0(R_0)$ into $\mathfrak{C}^1(R_0)$, and
from $\mathfrak{C}^{0,\alpha}(R_0)$ into
$\mathfrak{C}^{2,\alpha}(R_0)$, these operator norms being independent of $\mathsf{a}$. 

ii) If $g \in \mathfrak{C}^{0,\alpha}(R_0)$,
then $Q=\mathsf{a}^2\mathfrak{K}^{(n)}g$ satisfies the equation
$$\frac{\partial^2Q}{\partial\varpi^2}
+\frac{n-2}{\varpi}\frac{\partial Q}{\partial\varpi}+
\frac{\partial^2Q}{\partial z^2}+g=0$$
on $\mathfrak{D}(R)=\{ (\varpi, z) | r=\sqrt{\varpi^2+z^2}<R\}$.
\end{Proposition}

\subsection{Solutions of an equation of Poisson type with a negative potential}

Let us prove the following

\begin{Proposition}
There is a bounded linear operator $\mathfrak{L}$ in $\mathfrak{C}^0(R_0)$ which enjoys the following properties:

i) $\mathfrak{L}$ is continuous from
$\mathfrak{C}^0(R_0)$ into $\mathfrak{C}^1(R_0)$, and
from $\mathfrak{C}^{0,\alpha}(R_0)$ into
$\mathfrak{C}^{2,\alpha}(R_0)$, the operator norms being independent of
$\mathsf{a}$.

ii) For $g\in \mathfrak{C}^{0,\alpha}(R_0)$, the function
$Q=\mathsf{a}^2\mathfrak{L}g$ satisfies
the equation
$$\frac{\partial^2Q}{\partial\varpi^2}+
\frac{1}{\varpi}\frac{\partial Q}{\partial\varpi}+
\frac{\partial^2Q}{\partial z^2}+\frac{4\pi\mathsf{G}}{\gamma-1}
\frac{\rho_{\mathsf{N}}}{u_{\mathsf{N}}}Q+g =0$$
in $\mathfrak{D}(R)$ and $Q(0,0)=0$.
\end{Proposition}

Proof. First we note that
$$4\pi\mathsf{G}\frac{\rho_{\mathsf{N}}}{u_{\mathsf{N}}}=
\mathsf{a}^{-2}(\Theta\vee 0)^{\frac{1}{\gamma-1}-1}. $$
Here $\Theta$ is the distorted Lane-Emden function
$\Theta(r/\mathsf{a}, \zeta; \frac{1}{\gamma-1}, \mathsf{b})$. Thus the equation to be solved for $W=Q^{\flat(3)}, G=\mathsf{a}^2g^{\flat(3)}$ is
$$\triangle W+\frac{1}{\gamma-1}
(\Theta^{\flat}\vee 0)^{\frac{1}{\gamma-1}-1}W+G=0 $$
and we require $W(O)=0$. This problem reduces to the integral equation
$$W=\mathcal{K}\Big[\frac{1}{\gamma-1}
(\Theta^{\flat}\vee 0)^{\frac{1}{\gamma-1}-1}W
+G\Big],$$
where
$$\mathcal{K}f=\mathcal{K}^{(3)}f-
(\mathcal{K}^{(3)}f)(O).$$
In other words, we define $\mathfrak{K}$ by
$$
(\mathfrak{K} F)^{\flat(3)}=\mathcal{K}F^{\flat(3)} =\mathcal{K}^{(3)}F^{\flat(3)}-\mathcal{K}^{(3)}F^{\flat(3)}(O).
$$
The integral equation reads
$$Q=\mathfrak{K}\Big[
\frac{1}{\gamma-1}(\Theta\vee 0)^{\frac{1}{\gamma-1}-1}Q+\mathsf{a}^2g\Big].$$
Since the operator
$$
\mathfrak{T}=I-\mathfrak{K}\Big[\frac{1}{\gamma-1}
(\Theta\vee0)^{\frac{1}{\gamma-1}-1}\cdot\Big]: 
 Q\mapsto Q-\mathfrak{K}\Big[\frac{1}{\gamma-1}
(\Theta\vee0)^{\frac{1}{\gamma-1}-1}Q\Big]
$$
maps $\mathfrak{C}^0(R_0)$ into
$\mathfrak{C}^1(R_0)$, it is a compact operator in $\mathfrak{C}^0(R_0)$. But 
,  as noted in Introduction, it is known that the kernel of $\mathfrak{T}$ is 
$\{0\}$, that is, 
$$h=
\mathfrak{K}\Big[\frac{1}{\gamma-1}
(\Theta\vee0)^{\frac{1}{\gamma-1}-1}h\Big], \quad h\in\mathfrak{C}^0(R_0)\quad
\Rightarrow h=0,$$
provided that $\mathsf{b}\leq \varepsilon_0$ is sufficiently small. 
See \cite{JJTM} and \cite{JJTM.Ext}. Therefore the inverse 
$\mathfrak{T}^{-1}$ turns out to be a bounded linear
operator in $\mathfrak{C}^0(R_0)$. Then $\mathfrak{L}=\mathfrak{T}^{-1}
\mathfrak{K}$ is the required operator. Actually,
if $g \in \mathfrak{C}^0$ and $Q=\mathsf{a}^2\mathfrak{L}g$, then
$$Q=
\mathfrak{K}\Big[\frac{1}{\gamma-1}
(\Theta\vee0)^{\frac{1}{\gamma-1}-1}Q+\mathsf{a}^2g\Big]
$$
holds so that $Q \in \mathfrak{C}^1$. Moreover, if
$g\in \mathfrak{C}^{0,\alpha}$, then 
$\displaystyle \frac{1}{\gamma-1}
(\Theta\vee 0)^{\frac{1}{\gamma-1}-1}Q\in \mathfrak{C}^{0,\alpha}$
implies $Q \in \mathfrak{C}^{2,\alpha}$, and $Q$ satisfies the equation
$$\mathsf{a}^2\Big[\frac{\partial^2}{\partial\varpi^2}
+\frac{1}{\varpi}\frac{\partial}{\partial\varpi}+
\frac{\partial^2}{\partial z^2}\Big]Q+
\frac{1}{\gamma-1}(\Theta\vee 0)^{\frac{1}{\gamma-1}-1}Q
+\mathsf{a}^2g=0.
$$
Clearly
\begin{align*}
&\|Q; \mathfrak{C}^1\|\leq C\mathsf{a}^2\|g;\mathfrak{C}^0\|, \\
&\|Q; \mathfrak{C}^{2,\alpha}\|\leq C\mathsf{a}^2\|g; \mathfrak{C}^{0,\alpha}\|
\end{align*}
with a constant $C$ independent of $\mathsf{a}$. $\square$.

\section{Existence of solutions}

We are going to find solutions $w, Y, X, V$ of \eqref{Eq.w}\eqref{Eq.Y}\eqref{Eq.X}\eqref{Eq.K1}{\eqref{Eq.K3}.\\

We have fixed a small positive number $\varepsilon_0$ and suppose that $\mathsf{b}\leq \varepsilon_0$.

Let us fix $R$  such that
\begin{equation}
2\xi_1\Big(\frac{1}{\gamma-1}\Big) \leq \Xi=\frac{R}{\mathsf{a}}<\Xi_0=\frac{R_0}{\mathsf{a}}.
\end{equation}
We suppose that with a fixed finite constant $C_0 (\geq 1)$ it holds that
\begin{equation}
u_{\mathsf{O}}\leq C_0, \quad \frac{1}{\mathsf{c}^2}\leq C_0, \quad \frac{1}{C_0}\leq \mathsf{a}. 
\end{equation}
Moreover, we fix a positive number $\delta_0$ such that
\begin{equation}
\frac{|w|}{\mathsf{c}^2} \leq u_{\mathsf{O}}\delta_0 \Rightarrow u_{\mathsf{N}}+
\frac{w}{\mathsf{c}^2 }< 0 
\quad\mbox{for}\quad 2\mathsf{a}\xi_1 \leq r \leq R_0.
\end{equation}
Actually it is sufficient to take $\delta_0 <-\max_{|\zeta|\leq 1}\Theta(2\xi_1,\zeta)$. Then 
the support of $\displaystyle u=u_{\mathsf{N}}+\frac{w}{\mathsf{c}^2}$ is included in
$\mathfrak{D}(2\mathsf{a}\xi_1)$ as that of $u_{\mathsf{N}}$, provided that 
$\displaystyle \frac{|w|}{\mathsf{c}^2} \leq u_{\mathsf{O}}\delta_0$.

We are keeping in mind that
\begin{equation}
\mathsf{a}^{-2}=4\pi
\mathsf{G}\Big(\frac{\mathsf{A}\gamma}{\gamma-1}\Big)^{-\frac{1}{\gamma-1}}
u_{\mathsf{O}}^{\frac{1}{\gamma-1}-1},
\end{equation}
\begin{equation}
\Omega^2=
2\pi\mathsf{G}\mathsf{b}\Big(\frac{\gamma-1}{\mathsf{A}\gamma}\Big)^{\frac{1}{\gamma-1}}
u_{\mathsf{O}}^{\frac{1}{\gamma-1}}
=C\mathsf{b}\mathsf{a}^{-2}u_{\mathsf{O}}\leq C\mathsf{a}^{-2}u_{\mathsf{O}},\label{5.5}
\end{equation}
and
\begin{align}
& \|u_{\mathsf{N}}; \mathfrak{C}^{2,\alpha}(R_0)\|\leq C u_{\mathsf{O}}, \quad
\|\rho_{\mathsf{N}}; \mathfrak{C}^{2,\alpha}(R_0)\|\leq C\mathsf{A}^{-\frac{1}{\gamma-1}}
u_{\mathsf{O}}^{\frac{1}{\gamma-1}}
\leq C'\mathsf{G}^{-1}\mathsf{a}^{-2}u_{\mathsf{O}}, \nonumber \\
&\|P_{\mathsf{N}}; \mathfrak{C}^{2,\alpha}(R_0)\|\leq C\mathsf{A}^{-\frac{1}{\gamma-1}}
u_{\mathsf{O}}^{\frac{\gamma}{\gamma-1}}
\leq C'\mathsf{G}^{-1}\mathsf{a}^{-2}u_{\mathsf{O}}^2, \nonumber \\
&\|\Phi_{\mathsf{N}}; \mathfrak{C}^{2,\alpha}(R_0)\|\leq C u_{\mathsf{O}}, \quad \|\Phi_{\mathsf{N}}'; \mathfrak{C}^{2,\alpha}(R_0)\|\leq Cu_{\mathsf{O}}.
\label{5.6}
\end{align}
Here  $C, C'$ stand for constants which depend upon $\gamma, \alpha$ only.\\

First, supposing that $V \in \mathfrak{C}^{0,\alpha}(R_0)$ is given,  we are going to solve 
the equations \eqref{Eq.w}\eqref{Eq.Y}\eqref{Eq.X} for unknown $w, Y, X$ by solving the integral equations
\begin{subequations}
\begin{align}
&w=\mathfrak{L}(g_a+2\Omega^2 Y_1+R_a), \label{5.7a}\\
&Y=\mathfrak{K}^{(5)}(g_b+R_b), \label{5.7b}\\
&X=\mathfrak{K}^{(4)}(g_c+R_c), \label{5.7c}
\end{align}
\end{subequations}
where
\begin{subequations}
\begin{align}
&g_a=-8\Phi_{\mathsf{N}}'\Omega^2+
4\pi\mathsf{G}(-2\Phi_{\mathsf{N}}u_{\mathsf{N}}+
\lambda_1\rho_{\mathsf{N}}u_{\mathsf{N}}+3P), \\
&g_b=\frac{8}{\varpi}\frac{\partial\Phi_{\mathsf{N}}'}{\partial\varpi}, \\
&g_c=-16\pi\mathsf{G} P_{\mathsf{N}}.
\end{align}
\end{subequations}

From \eqref{5.6}, we see
\begin{align*}
&\|g_a;\mathfrak{C}^{0,\alpha}(R_0)\|\leq C\mathsf{a}^{-2}u_{\mathsf{O}}^2, \\
&\|g_b; \mathfrak{C}^{0,\alpha}(R_0)\|\leq C\mathsf{a}^{-2}u_{\mathsf{O}}, \\
&\|g_c; \mathfrak{C}^{0,\alpha}(R_0)\|\leq C\mathsf{a}^{-2}u_{\mathsf{O}},
\end{align*}
so that
\begin{align*}
&\|\mathfrak{L}g_a; \mathfrak{C}^{2,\alpha}(R_0)\|\leq C u_{\mathsf{O}}^2, \\
&\|\mathfrak{K}^{(5)}g_b; \mathfrak{C}^{2,\alpha}(R_0)\|\leq Cu_{\mathsf{O}}, \\
&\|\mathfrak{K}^{(4)}g_c; \mathfrak{C}^{2,\alpha}(R_0)\|\leq C u_{\mathsf{O}}^2.
\end{align*}\\

Given $w, Y, X \in \mathfrak{C}^{0,\alpha}(R_0)$, we evaluate $g_a, R_a, R_b, R_c$ by them, and we put
\begin{subequations}
\begin{align}
&\tilde{w}=\mathfrak{L}(g_a+2\Omega^2\tilde{Y}_1+R_a), \label{5.9a}\\
&\tilde{Y}=\mathfrak{K}^{(5)}(g_b+R_b), \label{5.9b}\\
&\tilde{X}=\mathfrak{K}^{(4)}(g_c+R_c).
\end{align}
\end{subequations}
Here $\tilde{Y}_1$ in the right hand side of \eqref{5.9a} means
$\displaystyle \varpi\frac{\partial\tilde{Y}}{\partial\varpi}+2\tilde{Y}$ given by $\tilde{Y}$
determined by \eqref{5.9b}.

Then our task is to find a fixed set of functions of the mapping
$(w,Y, X) \mapsto (\tilde{w}, \tilde{Y},\tilde{X})$.\\

In order to estimate $R_a$, we shall use the following estimates for 
 the function $H_{\rho}(w)$ is defined by
\begin{equation}
H_{\rho}(w)=f_{\rho}\Big(u_{\mathsf{N}}+\frac{w}{\mathsf{c}^2}\Big)-f_{\rho}(u_{\mathsf{N}})-Df_{\rho}(u_{\mathsf{N}})\frac{w}{\mathsf{c}^2}
\end{equation}
from the function $f_{\rho}(u)$, which gives $\rho_{\mathsf{N}}$ from
$u_{\mathsf{N}}$, 
\begin{equation}
f_{\rho}(u)=
\Big(\frac{\gamma-1}{\mathsf{A}\gamma}\Big)^{\frac{1}{\gamma-1}}
(u\vee 0)^{\frac{1}{\gamma-1}}:
\end{equation}

\begin{Proposition}\label{Proposition11}
i)  If $w \in \mathfrak{C}^0(R_0)$, then
\begin{equation}
\|H_{\rho}(w); \mathfrak{C}^{0}(R_0)\|\leq C
\Big(\frac{\|w;\mathfrak{C}^0\|}{\mathsf{c}^2}\Big)^{\frac{1}{\gamma-1}},\label{5.12}
\end{equation}
and
\begin{equation}
\mathsf{G}\mathsf{c}^2
\|H_{\rho}(w); \mathfrak{C}^{0}(R_0)\|\leq C\mathsf{a}^{-2}
\Big(\frac{1}{\mathsf{c}^2}\Big)^{\frac{1}{\gamma-1}-1}
u_{\mathsf{O}}^{-\frac{1}{\gamma-1}+1}
\|w;\mathfrak{C}^0\|^{\frac{1}{\gamma-1}}.\label{5.13}
\end{equation}

ii)  If $w_1,w_2\in \mathfrak{C}^0(R_0)$, then
\begin{equation}
\|H_{\rho}(w_2)-H_{\rho}(w_1); \mathfrak{C}^0(R_0)\|
\leq C \Big(\frac{1}{\mathsf{c}^2}(\|w_1;\mathfrak{C}^0\|+
\|w_2;\mathfrak{C}^0\|)\Big)^{\frac{1}{\gamma-1}-1}
\cdot
\frac{\|w_2-w_1;\mathfrak{C}^0\|}{\mathsf{c}^2},\label{5.14}
\end{equation}
and
\begin{align}
\mathsf{G}\mathsf{c}^2
\|H_{\rho}(w_2)-H_{\rho}(w_1); \mathfrak{C}^0(R_0)\|
&\leq C \mathsf{a}^{-2}
\Big(\frac{1}{\mathsf{c}^2}\Big)^{\frac{1}{\gamma-1}-1}
u_{\mathsf{O}}^{-\frac{1}{\gamma-1}+1}
(\|w_1;\mathfrak{C}^0\|+
\|w_2;\mathfrak{C}^0\|)\Big)^{\frac{1}{\gamma-1}-1} \times \nonumber \\
&\times
\|w_2-w_1;\mathfrak{C}^0\|.\label{5.15}
\end{align}

iii)  If $w \in \mathfrak{C}^1(R_0)$, then
$H_{\rho}(w)\in \mathfrak{C}^{0, \alpha}(R_0)$ and
\begin{align}
\|H_{\rho}(w);\mathfrak{C}^{0,\alpha}(R_0)\|&\leq C
u_{\mathsf{O}}^{\frac{1}{\gamma-1}-1}\Big[
\Big(1+\frac{\|w;\mathfrak{C}^1\|}{\mathsf{c}^2u_{\mathsf{O}}}
\Big)^{\frac{1}{\gamma-1}-1}
\frac{\|w;\mathfrak{C}^0\|}{\mathsf{c}^2} + \nonumber \\
&+\frac{\|w;\mathfrak{C}^0\|}{\mathsf{c}^2u_{\mathsf{O}}}
\frac{\|w;\mathfrak{C}^{0,\alpha}\|}{\mathsf{c}^2}\Big],\label{5.16}
\end{align}
and
\begin{align}
\mathsf{G}\mathsf{c}^2
\|H_{\rho}(w);\mathfrak{C}^{0,\alpha}(R_0)\|&\leq C
\mathsf{a}^{-2}\Big[
\Big(1+\frac{\|w;\mathfrak{C}^1\|}{\mathsf{c}^2u_{\mathsf{O}}}
\Big)^{\frac{1}{\gamma-1}-1}
\|w;\mathfrak{C}^0\| + \nonumber \\
&+
\frac{1}{\mathsf{c}^2}\frac{\|w;\mathfrak{C}^0\|}{u_{\mathsf{O}}}
\|w;\mathfrak{C}^{0,\alpha}\|\Big].\label{5.17}
\end{align}

iv) If $w_1,w_2 \in \mathfrak{C}^1(R_0)$, then
\begin{align}
\|H_{\rho}(w_2)-H_{\rho}(w_1);\mathfrak{C}^{0,\alpha}\|&\leq C
u_{\mathsf{O}}^{\frac{1}{\gamma-1}-1}\Big[\Big(1+\frac{1}{\mathsf{c}^2u_{\mathsf{O}}}
(\|w_1;\mathfrak{C}^1\|+\|w_2;\mathfrak{C}^1\|\Big)^{\frac{1}{\gamma-1}-1}
\frac{\|w_2-w_1;\mathfrak{C}^0\|}{\mathsf{c}^2} + \nonumber \\
&+\frac{1}{\mathsf{c}^2u_{\mathsf{O}}}(\|w_1;\mathfrak{C}^0\|+\|w_2;\mathfrak{C}^0\|)
\frac{\|w_2-w_1;\mathfrak{C}^{0,\alpha}\|}{\mathsf{c}^2}\Big],
\end{align}
and
\begin{align}
\mathsf{G}\mathsf{c}^2\|H_{\rho}(w_2)-H_{\rho}(w_1);\mathfrak{C}^{0,\alpha}\|&\leq C
\mathsf{a}^{-2}\Big[\Big(1+\frac{1}{\mathsf{c}^2u_{\mathsf{O}}}
(\|w_1;\mathfrak{C}^1\|+\|w_2;\mathfrak{C}^1\|\Big)^{\frac{1}{\gamma-1}-1}
{\|w_2-w_1;\mathfrak{C}^0\|} + \nonumber \\
&+\frac{1}{\mathsf{c}^2u_{\mathsf{O}}}(\|w_1;\mathfrak{C}^0\|+\|w_2;\mathfrak{C}^0\|)
{\|w_2-w_1;\mathfrak{C}^{0,\alpha}\|}\Big]. \label{5.19}
\end{align}
\end{Proposition}

Proof is elementary, if we use the expression
$$
H_{\rho}(w)=\int_0^1\Big(Df_{\rho}\Big(u_{\mathsf{N}}+t\frac{w}{\mathsf{c}^2}\Big)-
Df_{\rho}(u_{\mathsf{N}})\Big)dt\cdot \frac{w}{\mathsf{c}^2}, $$
and
$$H_{\rho}(w_2)-H_{\rho}(w_1)=
\int_0^1
\Big(Df_{\rho}\Big(u_{\mathsf{N}}+\frac{1}{\mathsf{c}^2}(w_1+t(w_2-w_1))\Big)
-Df_{\rho}(u_{\mathsf{N}})\Big)dt\cdot\frac{w_2-w_1}{\mathsf{c}^2}.
$$
Here, of course,
$$
Df_{\rho}(u)=\frac{1}{\gamma-1}
\Big(\frac{\gamma-1}{\mathsf{A}\gamma}\Big)^{\frac{1}{\gamma-1}}
(u\vee 0)^{\frac{1}{\gamma-1}-1}.
$$\\

Let us fix $V \in \mathfrak{C}^{0,\alpha}(R_0)$ such that $\|V; \mathfrak{C}^{0,\alpha}(R_0)\|\leq u_{\mathsf{O}}^2M$ and suppose
\begin{equation}
\frac{u_{\mathsf{O}}M}{\mathsf{c}^2}\leq C_0.
\end{equation}\\

Let us consider $w, Y, X$ such that
\begin{equation}
\|w;\mathfrak{C}^1\|\leq u_{\mathsf{O}}^2B,\quad
\|Y;\mathfrak{C}^1\|\leq u_{\mathsf{O}}B,
\quad
\|X;\mathfrak{C}^1\|\leq u_{\mathsf{O}}^2B,
\label{X1}
\end{equation}
and
\begin{equation}
\|w;\mathfrak{C}^{2,\alpha}\|\leq u_{\mathsf{O}}^2B^*,\quad
\|Y;\mathfrak{C}^{2,\alpha}\|\leq u_{\mathsf{O}}B^*,
\quad
\|X;\mathfrak{C}^{2,\alpha}\|\leq u_{\mathsf{O}}^2B^*, \label{X2}
\end{equation}
with \begin{equation}
B\leq B^*.
\end{equation}
Suppose
\begin{equation}
\frac{u_{\mathsf{O}}B}{c^2}\leq\delta_0. \label{Bd0}
\end{equation}\\

Using \eqref{5.13}, we have
\begin{align}
\|\tilde{w};\mathfrak{C}^1\|&\leq C_1\Big[u_{\mathsf{O}}^2+
u_{\mathsf{O}}\|\tilde{Y}_1;\mathfrak{C}^1\|+\frac{u_{\mathsf{O}}^3}{\mathsf{c}^2}
(B(1+B+M+B^*)+M)+ \nonumber \\
&+\Big(\frac{u_{\mathsf{O}}B}{\mathsf{c}^2}\Big)^{\frac{1}{\gamma-1}-1}u_{\mathsf{O}}^2B\Big].
\end{align}
Here we have used the estimate
$$\Big\|\frac{1}{\varpi}\frac{\partial Q}{\partial\varpi}; \mathfrak{C}^0\Big\|\leq C
\mathsf{a}^{-2}\|Q;\mathfrak{C}^2\|, $$
which can be verified by the identity
$$\mathsf{a}^2\frac{1}{\varpi}\frac{\partial Q}{\partial\varpi}=
\frac{1}{n-2}\sum_{j=2}^{n-1}
\frac{\partial^2Q^{\flat}}{\partial\xi_j^2}\Big|_{\xi=(\varpi/\mathsf{a},0,\cdots )}.
$$

On the other hand, using \eqref{5.17}, we have
\begin{equation}
\|\tilde{w};\mathfrak{C}^{2,\alpha}\|\leq C_1
\Big[ u_{\mathsf{O}}^2+u_{\mathsf{O}}\|\tilde{Y}; \mathfrak{C}^{2,\alpha}\|+
\frac{u_{\mathsf{O}}^3}{\mathsf{c}^2}(B^*(1+B^*+M)+M)+
u_{\mathsf{O}}^2B\Big].
\end{equation}

Clearly we have
\begin{align}
&\|\tilde{Y};\mathfrak{C}^1\|\leq C_1\Big[u_{\mathsf{O}}+
\frac{u_{\mathsf{O}}^2}{\mathsf{c}^2}(B+B^*)\Big], \\
&\|\tilde{Y};\mathfrak{C}^{2,\alpha}\|\leq C_1\Big[
u_{\mathsf{O}}+\frac{u_{\mathsf{O}}^2}{\mathsf{c}^2}B^*\Big], \\
&\|\tilde{X};\mathfrak{C}^1\|\leq
\|\tilde{X};\mathfrak{C}^{2,\alpha}\|\leq C_1
\Big[u_{\mathsf{O}}^2+\frac{u_{\mathsf{O}}^3}{\mathsf{c}^2}(B+M)\Big].
\end{align}\\

Then we can find $B, B^*$ such that $B\leq B^*$ which satisfy
\eqref{Bd0} and
\begin{subequations}
\begin{align}
&C_1\Big[1+C_1\Big(1+\frac{u_{\mathsf{O}}}{\mathsf{c}^2}B^*)\Big) + \nonumber \\
&+\frac{u_{\mathsf{O}}}{\mathsf{c}^2}(B(1+B+M+B^*)+M)+
\Big(\frac{u_{\mathsf{O}}B}{\mathsf{c}^2}\Big)^{\frac{1}{\gamma-1}-1}
B\Big]\leq B, \\
&C_1\Big[1+C_1\Big(1+\frac{u_{\mathsf{O}}}{\mathsf{c}^2}(B+B^*)\Big) + \nonumber \\
&+\frac{u_{\mathsf{O}}}{\mathsf{c}^2}(B(1+B+M+B^*)+M)+
B\Big]\leq B^*, \\
&C_1\Big[1+\frac{u_{\mathsf{O}}}{\mathsf{c}^2}(B+B^*)\Big]\leq B, \\
&C_1\Big[1+\frac{u_{\mathsf{O}}}{\mathsf{c}^2}B^*\Big] \leq B^*, \\
&C_1\Big[1+\frac{u_{\mathsf{O}}}{\mathsf{c}^2}(B+M)\Big]\leq B\leq B^*,
\end{align}
\end{subequations}
provided that the following assumption {\bf (C)} holds: \\

{\bf (C): \hspace{8mm} $u_{\mathsf{O}}/\mathsf{c}^2$ is sufficiently small
, say, $u_{\mathsf{O}}/\mathsf{c}^2 \leq \delta_1$, $\delta_1 (\leq 1)$ being a small positive number.} \\

Actually we can take
$$B=2C_1(1+C_1),\quad B^*=2C_1B=4C_1^2(1+C_1). $$

Then  the mapping $(w,Y,X) \mapsto (\tilde{w}, \tilde{Y}, \tilde{X})$ maps  the space $\mathfrak{X}$:
\begin{equation}
\mathfrak{X}:=\{(w,Y,X)\in \mathfrak{C}^{2,\alpha}(R_0)
\times\mathfrak{C}^{2,\alpha}(R_0)
\times\mathfrak{C}^{2,\alpha}(R_0)
\  |\  
\eqref{X1}, \eqref{X2} \mbox{hold}. \}
\end{equation}
into itself.\\

\begin{Remark}\label{Remark3}
When $\displaystyle\frac{1}{\gamma-1}>2 (\Leftrightarrow \gamma<\frac{3}{2})$, then the above roundabout plan of argument using both $\|\cdot; \mathfrak{C}^1\|$ and
$\|\cdot; \mathfrak{C}^{2,\alpha}\|$ is not necessary,
 and a simple argument can work.  Actually, if  
$\displaystyle\frac{1}{\gamma-1}>2$, then 
we have the expression
$$H_{\rho}(w)=\int_0^1
(1-t)D^2f_{\rho}\Big(u_{\mathsf{N}}+t\frac{w}{\mathsf{c}^2}\Big)dt\cdot \frac{w}{\mathsf{c}^2},$$
while $$
D^2f_{\rho}(u)=\Big(\frac{\gamma-1}{\mathsf{A}\gamma}\Big)^{\frac{1}{\gamma-1}}
\frac{2-\gamma}{(\gamma-1)^2}(u\vee 0)^{\frac{1}{\gamma-1}-2}$$
is continuous of $u$, and we see
$$\mathsf{G}\mathsf{c}^2\|H_{\rho}(w);\mathfrak{C}^{0,\alpha'}\|
\leq C\mathsf{a}^{-2}\frac{1}{\mathsf{c}^2u_{\mathsf{O}}}\|w;\mathfrak{C}^{0,\alpha'}\|^2,$$
which simplifies Proposition~\ref{Proposition11}, iii), \eqref{5.17}, provided that
$$\frac{\|w;\mathfrak{C}^1\|}{\mathsf{c}^2u_{\mathsf{O}}}\lesssim 1.$$
Therefore a simple argument using only
$\|\cdot; \mathfrak{C}^{2,\alpha'}\|$ can work. Here $\alpha'$ is taken as
$$0<\alpha' <\min\Big\{\frac{1}{\gamma-1}-2, 1\Big\}.$$
\end{Remark}

--\\

We are going to show that the mapping
$(w,Y,X)\mapsto (\tilde{w}, \tilde{Y}, \tilde{X})$ is a contraction with respect to a suitable norm.

Let us denote $U=(w,Y,X), \tilde{U}=(\tilde{w},\tilde{Y},\tilde{X})$ and so on. We consider the norms
\begin{subequations}
\begin{align}
&\mathcal{N}(U):=\max\{\|w;\mathfrak{C}^1\|, \|u_{\mathsf{O}}Y;\mathfrak{C}^1\|, \|X;\mathfrak{C}^1\|\}, \\
&\mathcal{N}^*(U):=\max\{\|w;\mathfrak{C}^{2,\alpha}\|, \|u_{\mathsf{O}}Y;\mathfrak{C}^{2,\alpha}\|, 
\|X;\mathfrak{C}^{2,\alpha}\|\}, \\
&\mathfrak{N}(U):=\mathcal{N}(U)+\kappa \mathcal{N}^*(U).
\end{align}
\end{subequations}
Here $\kappa$ is a positive constant specified later.

Using \eqref{5.15}, we see
\begin{equation}
\mathcal{N}(\tilde{U}_2-\tilde{U}_1)\leq C_2
\Big[\Big(\frac{u_{\mathsf{O}}}{\mathsf{c}^2}\Big)^{\alpha}\mathcal{N}(U_2-U_1)+
\frac{u_{\mathsf{O}}}{\mathsf{c}^2}\mathcal{N}^*(U_2-U_1)\Big].
\end{equation}

Using \eqref{5.19}, we see
\begin{equation}
\mathcal{N}^*(\tilde{U}_2-\tilde{U}_1)
\leq C_2\Big[\mathcal{N}(U_2-N_1)+\frac{u_{\mathsf{O}}}{\mathsf{c}^2}\mathcal{N}^*(U_2-U_1)\Big].
\end{equation}

Therefore we see
\begin{equation}
\mathfrak{N}(\tilde{U}_2-\tilde{U}_1)\leq K
[\mathcal{N}(U_2-U_1)+\kappa' \mathcal{N}^*(U_2-U_1)],
\end{equation}
where
\begin{align}
&K=C_2\Big[\Big(\frac{u_{\mathsf{O}}}{\mathsf{c}^2}\Big)^{\alpha}+\kappa\Big], \\
&\kappa'=
{\frac{u_{\mathsf{O}}}{\mathsf{c}^2}(1+\kappa)}\Big[\Big(\frac{u_{\mathsf{O}}}{\mathsf{c}^2}\Big)^{\alpha}+\kappa\Big]^{-1}.
\end{align}
Let us take
\begin{equation}
\kappa=2\Big(\frac{u_{\mathsf{O}}}{\mathsf{c}^2}\Big)^{1-\alpha},
\end{equation}
provided that
\begin{equation}
2\Big(\frac{u_{\mathsf{O}}}{\mathsf{c}^2}\Big)^{1-\alpha}
\leq 2\delta_1^{1-\alpha} \leq 1.
\end{equation}
Then we have $\kappa' \leq \kappa$ so that
\begin{equation}
\mathfrak{N}(\tilde{U}_2-\tilde{U}_1)\leq K\mathfrak{N}(U_2-U_1).
\end{equation}
But, under the assumption {\bf (C) } with sufficiently small $\delta_1$, we have 
$$K=
C_2\Big[\Big(\frac{u_{\mathsf{O}}}{\mathsf{c}^2}\Big)^{\alpha}+
2\Big(\frac{u_{\mathsf{O}}}{\mathsf{c}^2}\Big)^{1-\alpha}     
\Big]
\leq C_2\Big[\delta_1^{\alpha}+2\delta_1^{1-\alpha}\Big] <1,$$
 say, $U\mapsto \tilde{U}$ is a contraction with respect to the norm
$\mathfrak{N}$. Note that  $\mathfrak{N}$ is equivalent to $\mathcal{N}^*$ since $\kappa \mathcal{N}^*\leq\mathfrak{N}\leq (1+\kappa)\mathcal{N}^*$.\\

\begin{Remark}
When $\displaystyle\frac{1}{\gamma-1}>2$ as Remark~\ref{Remark3}, the argument can be simplified by using the expression
\begin{align*}
&H_{\rho}(w_2)-H_{\rho}(w_1)= \\
&=\int_0^1\int_0^1
D^2f_{\rho}\Big(u_{\mathsf{N}}+\frac{s}{\mathsf{c}^2}((1-t)w_1+tw_2)\Big)ds\frac{(1-t)w_1+tw_2}{\mathsf{c}^2}dt\cdot
\frac{w_2-w_1}{\mathsf{c}^2}.
\end{align*}
Then Proposition~\ref{Proposition11}, iv), \eqref{5.19} is simplified as
\begin{align*}
&\mathsf{G}\mathsf{c}^2\|H_{\rho}(w_2)-H_{\rho}(w_1);
\mathfrak{C}^{0,\alpha'}\| \leq \\
&\leq C\mathsf{a}^{-2}
\frac{1}{\mathsf{c}^2u_{\mathsf{O}}}
(\|w_1;\mathfrak{C}^0\|+\|w_2,\mathfrak{C}^0\|)
\|w_2-w_1; \mathfrak{C}^{0,\alpha'}\|,
\end{align*}
provided that
$$
\frac{1}{\mathsf{c}^2u_{\mathsf{O}}}
(\|w_1;\mathfrak{C}^1\|+\|w_2,\mathfrak{C}^1\|)\lesssim 1.$$
\end{Remark}

--\\

Summing up, we can claim

\begin{Theorem}\label{Th1}
There is a unique solution $(w,Y, X)$ of
\eqref{5.7a}\eqref{5.7b}\eqref{5.7c} in $\mathfrak{X}$ for any given
$V \in \mathfrak{C}^{0,\alpha}(R_0)$ with
$\|V;\mathfrak{C}^{0,\alpha}(R_0)\|\leq u_{\mathsf{O}}^2M$.
\end{Theorem}

\begin{Definition}
We shall denote the solution of Theorem~\ref{Th1} by
$U=(w,Y,X)=\mathcal{S}(V) $.
\end{Definition}

Note that $R_b$ does not depend on $V$ and the dependence of $R_a, R_c$ on $V$ occurs only through
$$\mathsf{c}^2(K'-F')=-\Phi_{\mathsf{N}}+\frac{1}{\mathsf{c}^2}(V+w).$$ 
Therefore, using both $\mathcal{N}$ and $\mathcal{N}^*$ to combine the estimates
of the form
\begin{align*}
\mathcal{N}(U_2-U_1)\leq & C\Big[\Big(\frac{u_{\mathsf{O}}}{\mathsf{c}^2}\Big)^{\alpha}\mathcal{N}(U_2-U_1)+
\frac{u_{\mathsf{O}}}{\mathsf{c}^2}\mathcal{N}^*(U_2-U_1) + \\
& + \frac{u_{\mathsf{O}}}{\mathsf{c}^2}\|V_2-V_1; \mathfrak{C}^0\|\Big], \\
\mathcal{N}^*(U_2-U_1)\leq & C\Big[
\frac{u_{\mathsf{O}}}{\mathsf{c}^2}\mathcal{N}^*(U_2-U_1) +
\mathcal{N}(U_2-U_1)+ \\
&+\frac{u_{\mathsf{O}}}{\mathsf{c}^2}\|V_2-V_1;\mathfrak{C}^{0,\alpha}\|\Big],
\end{align*}
where $U_{\mu}=\mathcal{S}(V_{\mu}), \mu=1,2$,
 it is easy to see

\begin{Theorem}\label{Th2}
There is a constant $C$ such that it holds for $V_{\mu} \in \mathfrak{C}^{0,\alpha}(R_0) $
with $\|V_{\mu}; \mathfrak{C}^{0,\alpha}(R_0)\|\leq u_{\mathsf{O}}^2M$, $\mu=1,2$, that
\begin{equation}
\mathcal{N}^*(\mathcal{S}(V_2)-\mathcal{S}(V_1))
\leq C\frac{u_{\mathsf{O}}}{\mathsf{c}^2}\|V_2-V_1;\mathfrak{C}^{0,\alpha}(R_0)\|,
\end{equation}
provided that $u_{\mathsf{O}}/\mathsf{c}^2 \leq \delta_1$ is sufficiently small.
\end{Theorem}

As the second step, we are going to solve the equations
\eqref{Eq.K1},\eqref{Eq.K3}. More precisely speaking, we are looking for 
$V\in \mathfrak{C}^{0,\alpha}(R_0)$ with
$\|V; \mathfrak{C}^{0,\alpha}(R_0)\|\leq u_{\mathsf{O}}^2 M$ such that $V$,
$U=\mathcal{S}(V)$
satisfy
\eqref{Eq.K1}, \eqref{Eq.K3}. 

Note that $V$ doe not appear in the right-hand side of \eqref{Eq.K1}, \eqref{Eq.K3}.

Let us denote $\mathcal{T}_1(V), \mathcal{T}_3(V)$
the right-hand sides of \eqref{Eq.K1}, \eqref{Eq.K3}, respectively,
evaluated by $U=\mathcal{S}(V)$, and let us put 
$\mathcal{T}(V)=\tilde{V}$ defined by
\begin{equation}
\tilde{V}:=
\int_0^z\mathcal{T}_3(V)(0,z')dz'+
\int_0^{\varpi}\mathcal{T}_1(V)(\varpi',z)d\varpi'.
\end{equation}
Thanks to Lemma~\ref{Lemma2}, if $V$ is a fixed point of the mapping
$\mathcal{T}$, it is the required solution
on $\mathfrak{D}(R)$.\\

We see clearly that there is a constant $C_3$ independent of $M$ such that
\begin{equation}
\|\mathcal{T}_j(V);\mathfrak{C}^{0,\alpha}\|\leq C_3\mathsf{a}^{-1}u_{\mathsf{O}}^2,\quad j=1,3,
\end{equation}
for $\|V;\mathfrak{C}^{0,\alpha}(R_0)'\|\leq u_{\mathsf{O}}^2M$ and $u_{\mathsf{O}}/\mathsf{c}^2\leq\delta_1 \ll 1$.  
    Then we have
\begin{equation}
\|\mathcal{T}(V); \mathfrak{C}^{0,\alpha}\|\leq  C_3u_{\mathsf{O}}^2. 
\end{equation}

So, taking $M$ large so that
\begin{equation}
C_3\leq M,\label{ZZ}
\end{equation}
we claim that under the mapping
$\mathcal{T}:V \mapsto \tilde{V}$ the functional set
\begin{equation}
\mathfrak{V}:=\{ V\in \mathfrak{C}^{0,\alpha}(R_0)\quad |
\quad
\|V; \mathfrak{C}^{0,\alpha}(R_0)\|\leq u_{\mathsf{O}}^2M\}
\end{equation}
is stable. Moreover, thanks to Theorem~\ref{Th2}, $\mathcal{T}$ is a contraction
with respect to the norm $\|\cdot; \mathfrak{C}^{0,\alpha}(R_0)\|$, provided $u_{\mathsf{O}}/\mathsf{c}^2\leq \delta_1$ is sufficiently small. 

Looking at the right-hand sides of \eqref{Eq.K1},\eqref{Eq.K3}, which are of class $\mathfrak{C}^{0,\alpha}$, we see  $V \in \mathfrak{C}^{1,\alpha}$. Then it follows that 
$X\in \mathfrak{C}^{3, \alpha}$, since 
$$
X=\mathfrak{K}^{(4)}[-16\pi
\mathsf{G}P_{\mathsf{N}}+R_c] $$
with
$$
-16\pi\mathsf{G}P_{\mathsf{N}}+R_c=
-16\pi\mathsf{G}
e^{2(K'-F')}P\Big(1+\frac{X}{\mathsf{c}^4}\Big) \in \mathfrak{C}^{1,\alpha}.
$$
In fact $V\in \mathfrak{C}^{1,\alpha}$ implies $K'-F' \in \mathfrak{C}^{1,\alpha}$ in view of \eqref{PS}.
As result $V$ turns out to be of class $\mathfrak{C}^{2,\alpha}$ in view of \eqref{Eq.K1},\eqref{Eq.K3} again.
Thus we have
\begin{Theorem}\label{Th3}
There is a solution $V\in
\mathfrak{C}^{2,\alpha}(R_0)$ of $V=\mathcal{T}(V)$
together with
$U=\mathcal{S}(V)$, provided that $u_{\mathsf{O}}/\mathsf{c}^2\leq\delta_1 \ll 1$. 
This is the unique solution in $\mathfrak{V}$. The equations
\eqref{Eq.K1}\eqref{Eq.K3} are satisfied on $\mathfrak{D}(R)$.
\end{Theorem}

Let us recall
\begin{align*}
&F'=
-\frac{\Omega^2}{2\mathsf{c}^2}\varpi^2+\frac{\Phi_{\mathsf{N}}}{\mathsf{c}^2}-\frac{w}{\mathsf{c}^4},
 \qquad K'=
-\frac{\Omega^2}{2\mathsf{c}^2}\varpi^2+\frac{V}{\mathsf{c}^4}, \\
&A'=
-\frac{\Omega}{\mathsf{c}}\varpi\Big(1+\frac{Y}{\mathsf{c}^2}\Big), \qquad \Pi=\varpi\Big(1+\frac{X}{\mathsf{c}^4}\Big), \qquad u=u_{\mathsf{N}}+\frac{w}{\mathsf{c}^2}.
\end{align*}

Note that, if $u_{\mathsf{O}}/\mathsf{c}^2 \ll 1$, then the assumption {\bf (B1)} holds,
since
\begin{align*}
e^{2F'}\Big(1-\frac{\Omega}{\mathsf{c}}A'\Big)^2-
e^{-2F'}\frac{\Omega^2}{\mathsf{c}^2}\Pi^2&=1+\frac{2\Phi_{\mathsf{N}}}{\mathsf{c}^2}+O
\Big(\frac{u_{\mathsf{O}}^2}{\mathsf{c}^4}\Big) \\
&=1+O\Big(\frac{u_{\mathsf{O}}}{\mathsf{c}^2}\Big),
\end{align*}
and the assumption {\bf (B2)} holds, since
$\displaystyle \frac{\partial\Pi}{\partial\varpi}=1+O\Big(
\frac{u_{\mathsf{O}}^2}{\mathsf{c}^4}\Big)$.

Thus the proof of the main result is complete.


\section{Physical vacuum boundary}

In this section, just in order to make sure, we shall see that the shape of the interior domain
$\mathfrak{R}:=\{\rho >0\}=\{u>0\}$ is given by a continuous function $R(\zeta)$ of $|\zeta|\leq 1$ so that 
$\mathfrak{R}=\{ r< R(\zeta)\}$. 

Now recall
\begin{equation}
u_{\mathsf{N}}=u_{\mathsf{O}}\Theta\Big(\frac{r}{\mathsf{a}}, \zeta; \frac{1}{\gamma-1},\mathsf{b}\Big),
\end{equation}
where $\Theta$ is the distorted Lane-Emden function constructed by
\cite{JJTM} and \cite{JJTM.Ext}. As proved in \cite[\S 4]{JJTM.Ext}, we may assume that 
\begin{equation}
\frac{\partial\Theta}{\partial r} <0\quad\mbox{for}\quad 0< r\leq \Xi_0.
\end{equation}

Therefore $\{\Theta(r,\zeta)>0\}$ is given by
$\{ r<\Xi_1(\zeta)\}$ with a continuous function
$\Xi_1(\zeta)=\Xi_1(\zeta;\frac{1}{\gamma-1},\mathsf{b})$ of
$|\zeta|\leq 1$ such that 
\begin{subequations}
\begin{align}
&\Xi_1(-\zeta)=\Xi_1(\zeta), \\
&|\Xi_1(\zeta)-\xi_1|\leq C\mathsf{b}, \\
&\Big|\frac{d\Xi_1}{d\zeta}\Big|\leq C\frac{\mathsf{b}}{\sqrt{1-\zeta^2}},
\end{align}
\end{subequations}
In other words,
\begin{equation}
u_{\mathsf{N}}(\varpi, z)>0 \quad\Leftrightarrow \quad r<\mathsf{a}\Xi_1(\zeta).\label{6.9}
\end{equation}
It follows from  \eqref{6.9} that,
for each fixed $\zeta$, the equation
\begin{equation}
u=u_{\mathsf{N}}+\frac{w}{\mathsf{c}^2}=0
\end{equation}
admits a unique solution $r=R(\zeta)$ in a neighborhood of $r=\mathsf{a}\Xi_1(\zeta)$, provided that $u_{\mathsf{O}}/\mathsf{c}^2$ is sufficiently small, since we have
$$|w|+\mathsf{a}\Big|\frac{\partial w}{\partial r}\Big|\leq Cu_{\mathsf{O}}^2$$
thanks to $\|w; \mathfrak{C}^1(R_0)\|\leq u_{\mathsf{O}}^2B$. Moreover $\zeta\mapsto R(\zeta)$ is continuous and
\begin{subequations}
\begin{align}
&|R(\zeta)-\mathsf{a}\Xi_1(\zeta)|\leq C\frac{\mathsf{a}u_{\mathsf{O}}^2}{\mathsf{c}^2}, \\
&\sqrt{1-\zeta^2}\Big|\frac{dR}{d\zeta}\Big|\leq C\Big(\mathsf{b}+\frac{u_{\mathsf{O}}^2}{\mathsf{c}^2}\Big).
\end{align}
\end{subequations}
The proof is standard one of the implicit function theorem. In fact we see
$$\frac{dR}{d\zeta}=-\Big(\frac{\partial}{\partial r}
\Big[u_{\mathsf{N}}+\frac{w}{\mathsf{c}^2}\Big]\Big)^{-1}
\frac{\partial}{\partial\zeta}\Big[u_{\mathsf{N}}+\frac{w}{\mathsf{c}^2}\Big]\Big|_{r=R(\zeta)}$$
and
$$\frac{d\Xi_1}{d\zeta}=-\mathsf{a}
\Big(\frac{\partial u_{\mathsf{N}}}{\partial r}\Big)^{-1}
\frac{\partial u_{\mathsf{N}}}{\partial\zeta}\Big|_{r=\mathsf{a}\Xi_1(\zeta)}=
\frac{1}{\sqrt{1-\zeta^2}}O(\mathsf{b}).
$$

Moreover we can claim that
\begin{equation}
\frac{\partial u}{\partial r}=\frac{\partial}{\partial r}\Big(u_{\mathsf{N}}+\frac{w}{\mathsf{c}^2}\Big)<0
\quad\mbox{for}\quad 0<r\leq\mathsf{a}\Xi_0, |\zeta|\leq 1.\label{P5.12}
\end{equation}
This guarantees that $\{u (\varpi, z)>0\}$ coincides with $\{r<R(\zeta)\}$.

Note that $\displaystyle\frac{\partial u}{\partial r}$ means
$\displaystyle\Big(\frac{\partial}{\partial r}\Big)_{\zeta=\mbox{Const.}}u$, and
that
$$\zeta=\Big(=\frac{z}{r}\Big)=\mbox{Const.}
\qquad \Leftrightarrow \qquad
\frac{z}{\varpi}=\Big(=\frac{\zeta}{\sqrt{1-\zeta^2}}\Big)=\mbox{Const.}.
$$
Thus \eqref{P5.12} says that $u$ is monotone decreasing with respect to $r$ along each half line $z/\varpi=\mbox{Const.}$\\

Let us consider a boundary point $\mathrm{P} \in \partial\mathfrak{R}=
\{r=R(\zeta)\}$. Let $\vec{N}$ denote the unit outer normal vector at 
$\mathrm{P}$ to the interior region $\mathfrak{R}=\{r<R(\zeta)\}$. 
We claim that
\begin{equation}
-\infty <\frac{\partial u}{\partial\vec{N}}\Big|_{\mathrm{P}}<0,\label{P5.16}
\end{equation}
that is,
\begin{equation}
-\infty <\frac{\partial}{\partial\vec{N}}\Big(\frac{dP}{d\rho}\Big)\Big|_{\mathrm{P}}<0.
\end{equation}
This means that the boundary $\partial\mathfrak{R}$ is a so called `{\bf physical vacuum boundary}'. 

In fact we have
$$\frac{\partial u}{\partial\vec{N}}=
\frac{1}{\sqrt{1+\frac{1-\zeta^2}{r^2}\Big(\frac{dR}{d\zeta}\Big)^2}}
\Big(\frac{\partial u}{\partial r}
-\frac{1-\zeta^2}{r^2}\frac{dR}{d\zeta}
\frac{\partial u}{\partial \zeta}\Big).
$$
See \cite[\S 6.3]{JJTM}. But
\begin{align*}
\frac{\partial u}{\partial r}&=\frac{u_{\mathsf{O}}}{\mathsf{a}}\Big[
\frac{d\theta}{dr}\Big|_{r=\xi_1}
+O(\mathsf{b})+O\Big(\frac{u_{\mathsf{O}}}{\mathsf{c}^2}\Big) \Big]\\
&=\frac{u_{\mathsf{O}}}{\mathsf{a}}\Big[-\frac{\mu_1}{\xi_1^2}+O(\mathsf{b})+
O\Big(\frac{u_{\mathsf{O}}}{\mathsf{c}^2}\Big)\Big]
\end{align*}
and
$$\sqrt{1-\zeta^2}\frac{\partial u}{\partial\zeta}=
\sqrt{1-\zeta^2}\frac{\partial u_{\mathsf{N}}}{\partial\zeta}+O\Big(\frac{u_{\mathsf{O}}^2}{\mathsf{c}^2}\Big)=u_{\mathsf{O}}\Big[
O(\mathsf{b})+O\Big(\frac{u_{\mathsf{O}}}{\mathsf{c}^2}\Big)\Big].
$$
Hence \eqref{P5.16} holds provided that $\mathsf{b}$ and $u_{\mathsf{O}}/\mathsf{c}^2$ are sufficiently small.

\section{Discussion}

We have constructed an axially symmetric metric on the domain $\mathfrak{D}=\{ r< R\}$, $R$ being arbitrarily large,
such that the support of the density
$\mathfrak{R}\cup\partial\mathfrak{R}=
\{ r \leq R(\zeta)\}$ is a compact subset of $\mathfrak{D}$. However we have not yet clarified what will happen when we continue the metric as long as possible
in the vacuum region for $r\geq R$. 

When we consider a  spherically symmetric metric, which is given by the 
Tolman-Oppenheimer-Volkoff equation, the interior metric can be patched with the Schwarzschild metric on the exterior vacuum region at the boundary $\partial\mathfrak{R}=\{ r=\mbox{Const.}\}$ in $C^2$-manner. See 
\cite[\textsf{SUPPLEMENTARY REMARK 4}]{KJM} and
\cite[Theorem 3]{CHTM}. But the author does not know an analogous statement to the Birkhoff's theorem on spherically symmetric metrics in the case of axially symmetric metrics. 

A discussion on this point will be given elsewhere. However let us note a remark as follows here.

  An exact formulation of the `matter-vacuum matching' problem
to find an {\bf asymptotically flat} vacuum exterior which can be patched to the interior solution can be found in the short review \cite{MacCallum} by three leading experts of this problem, and the papers \cite{Mars1999} and \cite{MacCallum2007}. Actually we note that, if we consider the equations \eqref{1.52a} \eqref{1.52b} \eqref{1.52c} on the vacuum region in which $\rho=0$, 
the equation \eqref{1.52c} admits the particular solution
$\Pi=\varpi$, and then the equations \eqref{1.52a}, \eqref{1.52b} 
are reduced to
\begin{align}
&\Big[\frac{\partial^2}{\partial\varpi^2}+\frac{1}{\varpi}\frac{\partial}{\partial\varpi}+
\frac{\partial^2}{\partial z^2}\Big]F'+
\frac{e^{4F'}}{\varpi}
\Big(\Big(\frac{\partial A'}{\partial \varpi}\Big)^2+
\Big(\frac{\partial A'}{\partial z}\Big)^2\Big)=0,\label{D1} \\
&\frac{\partial}{\partial \varpi}\Big(\frac{e^{4F'}}{\varpi}
\frac{\partial A'}{\partial \varpi}\Big)+\frac{\partial}{\partial z}
\Big(\frac{e^{4F'}}{\varpi}\frac{\partial A'}{\partial z}\Big)=0. \label{D2}
\end{align}
But it follows from \eqref{D2} that there exists a function $B'(\varpi, z)$
such that 
\begin{equation}
\frac{\partial B'}{\partial \varpi}=-\frac{e^{4F'}}{\varpi}
\frac{\partial A'}{\partial z},\quad
\frac{\partial B'}{\partial z}=\frac{e^{4F'}}{\varpi}
\frac{\partial A'}{\partial \varpi}. \label{D3}
\end{equation}
Then \eqref{D1} reads
\begin{equation}
\Big[\frac{\partial^2}{\partial\varpi^2}+\frac{1}{\varpi}\frac{\partial}{\partial\varpi}+
\frac{\partial^2}{\partial z^2}\Big]F'+
\frac{e^{-4F'}}{2}
\Big(\Big(\frac{\partial B'}{\partial \varpi}\Big)^2+
\Big(\frac{\partial B'}{\partial z}\Big)^2\Big)=0,\label{D4} 
\end{equation}
and the consistency condition for the existence of $A'$ reads
\begin{equation}
\Big[\frac{\partial^2}{\partial\varpi^2}+\frac{1}{\varpi}\frac{\partial}{\partial \varpi}
+\frac{\partial^2}{\partial z^2}\Big]B'
-4\Big(\frac{\partial F'}{\partial\varpi}\frac{\partial B'}{\partial\varpi}
+\frac{\partial F'}{\partial z}\frac{\partial B'}{\partial z}\Big)=0. \label{D5}
\end{equation}
So, the system of equations \eqref{D4}\eqref{D5} governs the `Ernst potential' $E'=e^{2F'}+\sqrt{-1}B'$. Let us consider $F$, and $B$ which is defined by \eqref{D3} by replacing $B', A'$ by $B, A$. Then as remarked in the footnote
of \cite[p. 9]{Meinel}, $F, B$ satisfy the same equations \eqref{D4}, \eqref{D5} in which $F', B'$ are replaced by
$F, B$. This system is nothing but 
\cite[(1)]{Mars1999} or \cite[(1)]{MacCallum}, \cite[(11)]{MacCallum2007}. ($U, \Omega$ there  should read $F, B$ here.) The Ernst potential $E=e^{2F}+\sqrt{-1}B$ obeys the `Ernst equation'
$$\mathfrak{Re} [E]\cdot
\Big[\frac{\partial^2}{\partial\varpi^2}+\frac{1}{\varpi}\frac{\partial}{\partial\varpi}
+\frac{\partial^2}{\partial z^2}\Big]E=
\Big(\frac{\partial E}{\partial\varpi}\Big)^2+
\Big(\frac{\partial E}{\partial z}\Big)^2,
$$
which is equivalent to the system \eqref{D4},\eqref{D5}. Then the asymptotic flatness condition reads
$$
F=1-\frac{\mathsf{G}M}{\mathsf{c}^2r}+O\Big(\frac{1}{r^2}\Big),
\quad
B=-\frac{2\mathsf{G}zJ}{\mathsf{c}^3r^3}+O\Big(\frac{1}{r^3}\Big)
\quad\mbox{as}\quad r\rightarrow +\infty
$$
with some constants $M, J$. Here, of course, $r=\sqrt{\varpi^2+z^2}$.
See \cite{MacCallum}.

Thus it is an important open problem to find an asymptotically flat vacuum metric which is patched to the metric constructed here on an bounded domain. 

But we should note that the metric of the form \eqref{2} with the co-rotating potential $F'$ cannot be global and asymptotically flat. Actually, if it is possible, it should hold at least that
$$F \rightarrow 0,\quad A\rightarrow 0, \quad K\rightarrow 0, \quad
\Pi \sim \varpi $$
as $ r\rightarrow +\infty$. See \cite[p. p.5, (1.12)]{Meinel}. Then 
the left-hand side of \eqref{B2}
$$
e^{2F}\Big(1+\frac{\Omega}{\mathsf{c}}A\Big)^2-
e^{-2F}\frac{\Omega^2}{\mathsf{c}^2}\Pi^2\quad \sim\quad 1-\frac{\Omega^2}{\mathsf{c}^2}\varpi^2
$$
turns out to be negative, provided that $\Omega \not=0$, as $\varpi \rightarrow \infty$,
and $F'$ cannot be real-valued for $e^{2F'} <0$. \\

However the formulation of the problem done in the work by U. Heilig \cite{Heilig} is quite different. In this work the basic equations comes from the  `frame theory' of J\"{u}rgen Ehlers
 (\cite{Ehlers}). This formulation of the Einstein's equations contains the perspective of a post-Newtonian approximation already. The basic exposition of this Ehlers' frame theory can be found in \cite{Lottermoser}, while an introductory review can be found in \cite{OliynykS}, and recent developments can be found in
\cite{Oliynyk2007}, \cite{Oliynyk2009}, \cite{Oliynyk2010}.

U. Heilig \cite{Heilig} claims that the metric is constructed globally, but the result seems not to be a solution of the matter-vacuum matching problem. For the metric
$ds^2=g_{\mu\nu}dx^{\mu}dx^{\nu}$ is claimed to satisfy
$$ g_{\mu\nu}-\eta_{\mu\nu} \in M_{2,\delta}^p(\mathbb{R}^3)\quad\mbox{with}\quad
p\geq 4, 0\leq \delta <-2+3\frac{p-1}{p}, $$
where $\eta_{\mu\nu}dx^{\mu}dx^{\nu}$ is the Minkowski metric, and the functional space 
$M_{2,\delta}^p(\mathbb{R}^3)$ 
consists of functions $f$ such that
$$ \sum_{|\ell|\leq 2}
\|(\sqrt{1+|x|^2})^{\delta+|\ell|}\partial^{\ell}f\|_{L^p(\mathbb{R}^3)} <\infty.$$
Therefore it seems that this does not guarantees that 
$$ g_{\mu\nu}-\eta_{\mu\nu}=O\Big(\frac{1}{|x|}\Big). $$

\vspace{20mm}

{\bf\Large Acknowledgment}\\

The author would like to express his sincere thanks to Professor Todd A. Oliynyk (Monash University) who kindly drew the attention of the author to the pioneering work by U. Heilig and related previous works. Sincere thanks of the author should be addressed also to Professor Marc Mars (Universidad de Salamanca), who kindly gave instructions on his works on the matter-vacuum matching 
problem during the Conference ``Between Geometry and Relativity'' held at Vienna, July 24-28, 2017.

This work has been done partially during the stay of the author at  the Erwin Schr\"{o}dinger 
International Institute for Mathematics and Physics at Vienna on July 24-29, 2017 in the program on ``Geometry and Relativity''. The author would like to  express his sincere thanks to ESI for the support and the hospitality during the stay,   
and to the Yamaguchi University Foundation for the financial support by a 
 research grant A1-2 (2017). The study was partially done during the stay of the author at Department of Mathematics of National University of Singapore on January 17- February 1, 2018. The author expresses his sincere thanks to Professor Shih-Hsien Yu 
for his kind invitation and discussions and to Department of Mathematics, NUS, for the financial support and hospitality during the stay.

The author would like to express his sincere thanks to the reviewer for taking pains of careful  reading of the manuscript containing many and complicated equations and for giving accurate suggestions of revisions.




\end{document}